%% file: main.tex
\documentclass{article}
\usepackage[a4paper, total={6in,9in}]{geometry}
\usepackage[utf8]{inputenc}
\usepackage[T1]{fontenc}
\usepackage{lmodern}

\usepackage[pdftex]{graphicx} 
\usepackage{subcaption}
\graphicspath{ {./images/} }
\usepackage{enumitem}
\usepackage[normalem]{ulem}
\usepackage[toc,page]{appendix}

\usepackage[bookmarks, bookmarksnumbered,
			colorlinks=true,
            linkcolor = blue,
            urlcolor  = blue,
            citecolor = teal]{hyperref}
\newcommand\fnsurl[1]{{\footnotesize\url{#1}}}

\usepackage{myquicksetup}
\mathtoolsset{showonlyrefs} 

\usepackage{mylivemacros}

\usepackage{l1interp} 

\usepackage{subfiles} 

\title{Tight bounds for minimum $\ell_1$-norm interpolation of noisy data}
\date{}

\usepackage[affil-it]{authblk}
\author[1]{Guillaume Wang\footnote{Equal contribution. Email addresses: konstantin.donhauser@ai.ethz.ch, guillaumewang21@gmail.com}}
\newcommand\CoAuthorMark{\footnotemark[\value{footnote}]} 
\author[1,2]{Konstantin Donhauser\protect\CoAuthorMark}
\author[1]{Fanny Yang}
\affil[1]{ETH Zurich, Department of Computer Science}
\affil[2]{ETH AI Center}

\hypersetup{
    pdftitle={Tight bounds for minimum l1-norm interpolation of noisy data},
    pdfsubject={math.ST; cs.IT; cs.LG; stat.ML},
    pdfauthor={Guillaume Wang, Konstantin Donhauser, Fanny Yang},
    pdfkeywords={}
}
\newif\ifarxiv
\arxivtrue
\newif\ifarx
\arxtrue
\begin{document}
\maketitle

\subfile{sections/abstract}

\subfile{sections/intro}

\subfile{sections/discussion}

\subfile{sections/proof_sketch}
\subfile{sections/proof_maintext}

\subfile{sections/conclusion}

\section*{Acknowledgments}

K.D.~is supported by the ETH AI Center and the ETH Foundations of Data Science. 
We would like to thank Geoffrey Chinot and Matthias L\"offler for insightful discussions.
Finally, we are very grateful to Sara van de Geer for her valuable advice. 
\bibliographystyle{alpha}
\bibliography{references}
\addcontentsline{toc}{section}{\refname} 

\newpage
\appendix


\subfile{sections_proof/l1_proof}

\end{document}

%% file: sections/abstract.tex
\begin{abstract}
We provide matching upper and lower bounds of order $\sigma^2/\log(d/n)$ for the prediction error of the minimum $\ell_1$-norm interpolator, a.k.a.\ basis pursuit. Our result is tight up to negligible terms when $d \gg n$, and is the first to imply asymptotic consistency of noisy minimum-norm interpolation for isotropic features and sparse ground truths. Our work complements the literature on ``benign overfitting'' for minimum $\ell_2$-norm interpolation, where asymptotic consistency can be achieved only when the features are effectively low-dimensional.

\end{abstract}

%% file: sections/intro.tex
\ifarx
\section{Introduction} 
\else
\section{INTRODUCTION} 
\fi
\label{sec:intro}

Recent experimental studies \cite{belkin_2019,zhang_2021} reveal that in the
modern high-dimensional regime, models that perfectly fit noisy
training data can still generalize well. The phenomenon stands in contrast to 
the classical wisdom that interpolating the data results in poor
statistical performance due to overfitting.  Many
theoretical papers have explored why, when, and to what extent
interpolation can be harmless for generalization, suggesting a coherent storyline:
High dimensionality itself can have a regularizing effect, in the
sense that it lowers the model's sensitivity to noise.
This intuition emerges from the fast-growing
literature studying min-$\ell_2$-norm interpolation in the
regression setting with input dimension $d$ substantially exceeding
sample size $n$ (see \cite{bartlett_2020,dobriban_2018} and
references therein).  Results and intuition for this setting also
extend to kernel methods \cite{ghorbani_2021,mei_2019}.  
\ifarx \else \\ \fi

However, a
closer look at this literature reveals that while high
dimensionality decreases the sensitivity to noise (error due to variance),
the prediction error generally does not vanish as $d,n\to\infty$. 
Indeed, the bottleneck for asymptotic consistency is a non-vanishing bias term which can only be avoided
when the features have low effective dimension
$\effdim = \trace{\Sigma}/\spnorm{\Sigma} \ll n$, where $\Sigma$ is the covariance matrix \cite{tsigler_2020}.
Therefore, current theory 
does not yet provide a convincing explanation
for why interpolating models generalize well for inherently high-dimensional input data. 
This work takes a step towards addressing this gap.

When the input data is effectively high-dimensional (\eg\ isotropic and $d \gg n$),
we generally cannot expect any 
data-driven estimator 
to generalize well
unless there is
underlying
structure that can be exploited.
In this paper, we hence focus on linear regression on isotropic Gaussian features with the simplest
structural assumption: sparsity of the ground truth in the
standard basis. For this setting, the $\ell_1$-penalized regressor (LASSO, \cite{tibshirani_1996})
achieves minimax optimal rates in the presence of noise
\cite{vandegeer_2008}, while basis pursuit (BP, \cite{chen_1998}) -- that is
min-$\ell_1$-norm interpolation -- generalizes well in
the noiseless case but is known to be very sensitive to noise
\cite{candes_2008,donoho_2006}.

Given recent results on high dimensionality decreasing sensitivity of
interpolators to noise, and classical results on the low bias of BP for
learning sparse signals, the following question naturally arises:
\begin{center}
    \emph{Can we consistently learn sparse ground truth functions 
    with minimum-norm interpolators on inherently %
    high-dimensional features?}
\end{center}
So far, upper bounds on 
the prediction error
of the BP
estimator
of the order of the noise level $O(\sigma^2)$
have been derived for isotropic Gaussian \cite{koehler_2021,ju_2020,wojtaszczyk_2010},
sub-exponential \cite{foucart_2014}, or heavy-tailed
\cite{chinot_2021,krahmer_2018} features.  
In the case of isotropic Gaussian features, even though
\ifarx the authors of the paper \cite{chinot_2021}
\else \citet{chinot_2021}
\fi
show a tight matching lower bound for adversarial noise, 
for \iid\ noise the best known results are not tight:
there is a gap between the non-vanishing upper bound $O(\sigma^2)$
\cite{wojtaszczyk_2010} and the lower bound $\Omega\left( \frac{\sigma^2}{\log(d/n)} \right)$ \cite{chatterji_2021,muthukumar_2020}.
For \iid\ noise, 
\ifarx the authors of the paper \cite{chinot_2021} 
\else \citet{chinot_2021}
\fi
conjecture that BP does not achieve consistency (see also \cite{koehler_2021}).

\paragraph{Contribution.}
We are the first to answer the above question in the affirmative.
Specifically, we show that for isotropic Gaussian features,
BP does in fact achieve asymptotic consistency
when $d$ grows superlinearly and subexponentially in~$n$,
disproving the recent conjecture by \cite{chinot_2021}.
Our result closes the aforementioned gap in the
literature on BP: We give matching upper and lower bounds of order~${ \frac{\sigma^2}{\log(d/n)}}$
on the prediction error
of the BP estimator, exact up to terms that are negligible when $d \gg n$.
Further, our proof technique is novel and may be of independent interest.

\paragraph{Structure of the paper.}
The rest of the article is structured as follows.  In
Section~\ref{sec:main_result}, we give our main result and
discuss its implications.
In Section~\ref{sec:proof_sketch}, we present a proof sketch
and provide insights on why our approach leads to tighter bounds than previous works.
We discuss the 
scope
of our assumptions and motivate future work in Section~\ref{subsec:future_work},
and conclude the paper in
Section~\ref{sec:conclusion}.


%% file: sections/discussion.tex
\ifarx
\section{Main result} 
\else
\section{MAIN RESULT} 

\fi
\label{sec:main_result}

In this section we state our main result, followed by a
discussion of its implications.
We consider a linear regression model with input vectors $x \in
\RR^d$ drawn from an isotropic Gaussian distribution $x \sim
\Normal(0,I_d)$, and response variable $y = \innerprod{\wgt}{x} +
\xi$, where $w^*$ is the ground truth to be estimated and $\xi \sim
\Normal(0,\sigma^2)$ is a noise term independent of $x$.
Given $n$ \iid\ random samples $(x_i,y_i)_{i=1}^n$, the goal is to estimate $\wgt$
and obtain a small prediction error
(or risk)
for the estimate $\what$
\begin{equation}
    \EE_{x,y} (\langle \what, x\rangle - y)^2  -\sigma^2
    = \| \what - \wgt\|_2^2
\end{equation}
where we subtract the irreducible error $\sigma^2$. Note that this is also exactly the $\ell_2$-error of the estimator.
We study the min-$\ell_1$-norm interpolator (or BP solution)
defined by 
\begin{equation}
    \what = \argmin_w \|w\|_1 
    ~~~\mathrm{such~that}~~~
    \forall i,~ \langle x_i, w \rangle = y_i.
\end{equation}


Our main result, Theorem~\ref{thm:main}, provides non-asymptotic matching upper and lower bounds for the prediction error of this estimator:

\begin{theorem} \label{thm:main} 
    Suppose 
    ${\norm{\wgt}_0 \leq \cwone \frac{n}{\log(d/n)^{5}}}$ 
    for some universal constant $\cwone>0$.
    There exist universal constants $\cn,\cl,\cd,\cone, \ctwo, \cprobgeoff>0$ such that,
    for any $n,d$ with
    $n \geq \cn$ and 
    $\cl n \log(n)^2 \leq d \leq \exp(\cd n^{1/5})$, the prediction error satisfies
    \begin{equation}
    \label{eq:mainbound}
        \abs{
            \| \what - \wgt\|_2^2 
             - \frac{\sigma^2}{\log(d/n)} 
        }
        \leq \cone \frac{\sigma^2}{\log(d/n)^{3/2}}
    \end{equation}
  with probability at least 
    $1- \ctwo \exp\left(-\frac{n}{\log(d/n)^5}\right)$  $- d \exp\left(-\cprobgeoff n \right)$
    over the draws of the dataset.
\end{theorem}

A proof sketch is presented in Section~\ref{sec:proof_sketch} and the full proof is given in \ifarx Section\else Appendix\fi~\ref{sec:proof_maintext}. We refer to Section~\ref{subsec:future_work} for a discussion on limitations of the assumptions. 

This theorem proves an exact
statistical rate with respect to the leading factor of order
$\frac{\sigma^2}{\log(d/n)}$ for the prediction error of the BP solution. 
Previous lower bounds
of order $\Omega \left(
\frac{\sigma^2}{\log(d/n)} \right)$ 
for the same distributional
setting (isotropic Gaussian features, \iid\ noise) only apply  under
more restrictive assumptions, such as the zero-signal case $\wgt = 0$
\cite{muthukumar_2020}, or assuming $d >n^4$
\cite{ju_2020}. 

On another note the best known  upper bounds are of constant order $O(\sigma^2)$
\cite{chinot_2021,wojtaszczyk_2010}. Our result both proves the lower bound in more generality and significantly improves the upper bound that matches the lower bound, showing that the lower bound is in fact tight.
An important implication of the upper bound is that BP achieves high-dimensional asymptotic consistency when $d = \omega(n)$,
thus disproving to a recent conjecture by \cite{chinot_2021}.


\paragraph{Dependency on $\wgt$.}
We note that the bound for the risk in Theorem~\ref{thm:main} is independent of the choice of $\wgt$ assuming that it is sparse (i.e, has bounded $\ell_0$-norm). Essentially, this arises from the well known fact that in the noiseless case ($\sigma =0$) we can achieve exact recovery \cite{candes_2008} of sparse ground truths. 
More generally, existing upper bounds for 
the prediction error of the BP estimator
for general ground truths $\wgt$ are of the form%
\footnote{The notation $a \lesssim b$ means that there exists a universal constant $c_1>0$ such that $  a \leq c_1 b$, and we write
$a \asymp b$ for $a \lesssim b$ and $b \lesssim a$.}
\begin{equation}
\label{eq:chinotbound}
    \norm{\what - \wgt}_2^2
    \lesssim 
        \frac{\|\xi\|_2^2}{n} + \|\wgt\|_1^2 \frac{\log(d/n)}{n}
\end{equation}
(see \eg\ \cite[Theorem~3.1]{chinot_2021}).
That is, they contain a first term reflecting the error due to overfitting of the noise $\xi$ which is independent of $\wgt$
and a second term which can be understood as the noiseless error only depending on $\wgt$ but not on the noise $\xi$. In fact, the authors of both papers show that assuming the ground truth is hard-sparse (bounded $\ell_0$-norm),  the second term on the RHS in Equation~\eqref{eq:chinotbound} can be avoided, resulting in the bound $ \norm{\what - \wgt}_2^2 \lesssim \frac{\|\xi\|_2^2}{n}$. Therefore, it is also not surprising that our tighter bound in Equation~\eqref{eq:mainbound} does not explicitly depend on $\wgt$.

\begin{figure*}[t]
    \centering
     \begin{subfigure}[b]{0.45\linewidth}
         \centering
    \includegraphics[width=\linewidth]{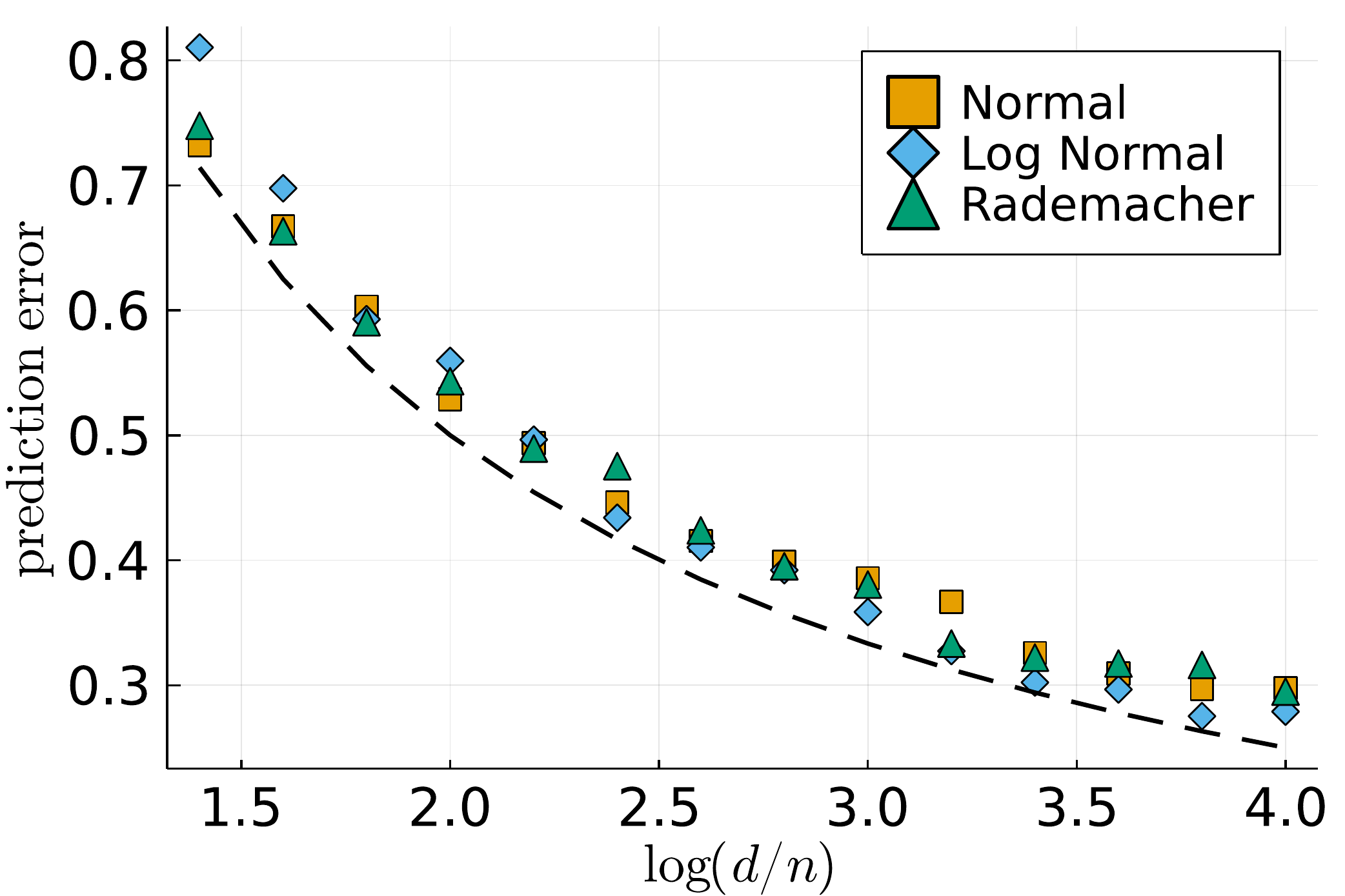}
         \caption{Prediction error of BP vs.\ $\log(d/n)$} \label{fig:distributions}
     \end{subfigure}
          \begin{subfigure}[b]{0.45\linewidth}
         \centering
        \includegraphics[width=\linewidth]{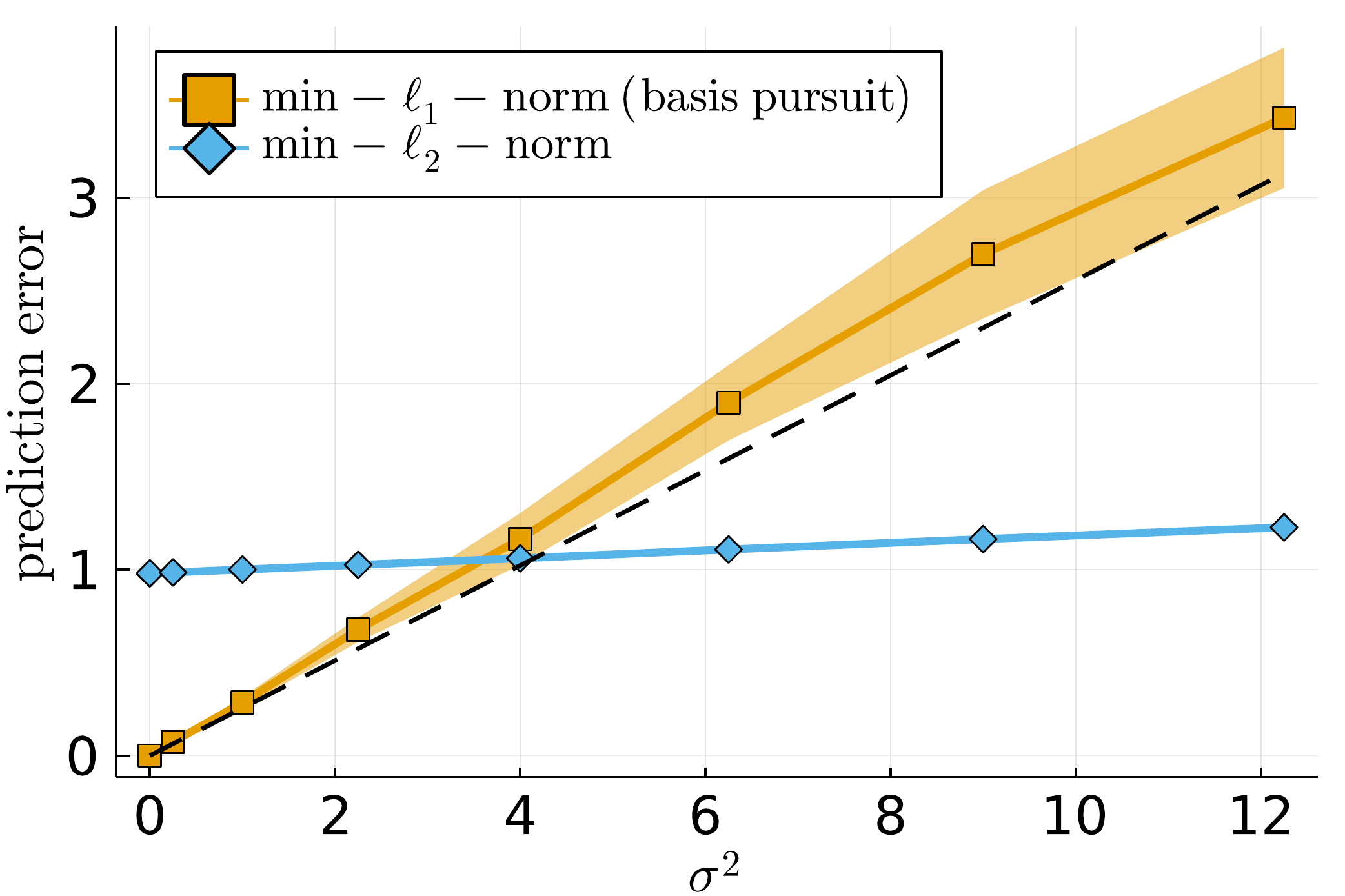}
         \caption{Prediction error vs.\ $\sigma^2$}
         \label{fig:tradeoff}
     \end{subfigure}
    \captionsetup{width=.92\linewidth}
    \caption{Prediction error as a function of
    (a) $\log(d/n)$ with varying $d$ and $n=400$ fixed, and 
    (b) $\sigma^2=1$ with $d=20000,n=400$.
    The features are generated by drawing from the isotropic 
    zero-mean and unit-variance
    (b) \emph{Normal} and (a)~\emph{Normal}, \emph{Log Normal} and \emph{Rademacher} distributions.
    For BP on Gaussian-distributed features (orange squares),
    the plots correctly reflect the theoretical rate  
    ${\frac{\sigma^2}{\log(d/n)}}$
    (dashed curve).
    See Section~\ref{subsec:numerical} for
    further details.}
    \label{fig:figure1}
    \vspace{-0.1in}
\end{figure*}

\subsection{Numerical simulations}
\label{subsec:numerical}

We now present numerical simulations illustrating Theorem~\ref{thm:main}. Figure~\ref{fig:distributions} shows the prediction error of BP plotted as a function of $\log(d/n)$ with varying $d$ and $n=400$ fixed, for isotropic inputs generated from the
zero-mean and unit-variance
\emph{Normal}, \emph{Log Normal} and \emph{Rademacher} distributions. For all three distributions, the prediction error closely follows the trend line
$\frac{\sigma^2}{\log(d/n)}$ (dashed curve). While Theorem~\ref{thm:main} only applies for Gaussian features, the figure suggests that this statistical rate of BP holds more generally (see discussion in Section~\ref{subsec:future_work}).

Figure~\ref{fig:tradeoff} shows the prediction error of the min-$\ell_1$-norm (BP) and min-$\ell_2$-norm interpolators as a function of the noise $\sigma^2$, for fixed $d$ and $n$. 
The prediction error of the former again aligns with the theoretical rate $\frac{\sigma^2}{\log(d/n)}$.
Furthermore, we observe that the min-$\ell_1$-norm
interpolator is sensitive to the noise level $\sigma^2$, while the min-$\ell_2$-norm interpolator has a similar (non-vanishing) prediction error across all values of $\sigma^2$.

For both plots we use $n=400$ and average the prediction error over 20 runs; in Figure~\ref{fig:tradeoff} we additionally show the standard deviation (shaded regions). The ground truth is $\wgt = (1,0,\cdots,0)$.
Finally, we choose $\sigma^2 =1$ in Figure~\ref{fig:distributions}, and $d = 20000$ in Figure~\ref{fig:tradeoff}.

\subsection{Implications and insights}
We now discuss 
further
high-level implications and insights that follow from Theorem~\ref{thm:main}.

\paragraph{High-dimensional asymptotic consistency.}
Our result proves consistency of BP for any asymptotic regime $d \asymp n^{\beta}$ with $\beta >1$.
In fact, we argue that those are the only regimes of interest.
For $d$ growing exponentially with $n$, known minimax lower bounds for sparse problems
of order $\frac{\sigma^2 s\log(d/s)}{n}$ (with $s \leq n$ the $\ell_0$-norm of the BP estimator), preclude consistency \cite{verzelen_2012}.
On the other hand, for linear growth ${d \asymp n}$, i.e., $\beta = 1$ -- studied in detail 
\ifarx in the paper \cite{li_2021} 
\else by \citet{li_2021}
\fi
--, the uniform prediction error lower bound $\frac{\sigma^2 n}{d-n}$ holding for all interpolators \cite{zhou20,muthukumar_2020} also forbids vanishing prediction error.
Note that for $d \asymp n^\beta$ ($\beta>1$), asymptotic consistency can also be achieved
by a carefully designed ``hybrid'' interpolating estimator \cite[Section 5.2]{muthukumar_2020}; contrary to BP, this estimator is not a minimum-norm interpolator, and is not structured (not $n$-sparse).

\paragraph{Trade-off between structural bias and  sensitivity to noise.}
As mentioned in the introduction,
our upper bound on the prediction error shows that, contrary to min-$\ell_2$-norm interpolation, BP is able to learn sparse signals 
in high dimensions
thanks to its structural bias
towards sparsity. 
However, our lower bound can be
seen as a tempering negative result: The prediction error decays only at a slow
rate of $\frac{\sigma^2}{\log(d/n)}$. 

Compared to min-$\ell_2$-norm interpolation, BP (min-$\ell_1$-norm interpolation) suffers from a higher sensitivity to noise, but possesses a more advantageous structural bias.
To compare the two methods' sensitivity to noise, consider the case $\wgt = 0$, where the prediction error purely reflects the effect of noise.
In this case, although both methods achieve vanishing error, the statistical rate for BP, $\frac{\sigma^2}{\log(d/n)}$, is much slower than that of min-$\ell_2$-norm interpolation, $\sigma^2 \max(\frac{1}{\sqrt{n}},\frac{n}{d})$ \cite[Theorem 3]{koehler_2021}.
Contrariwise, to compare the effect of structural bias, consider the noiseless case with a 
non-zero ground truth.
It is well
known that BP successfully learns sparse signals \cite{candes_2008},
while min-$\ell_2$-norm interpolation always fails to learn the ground truth due to the
lack of any corresponding structural bias.

Thus, there appears to be a trade-off between structural bias and sensitivity to noise: BP benefits from a strong
structural bias, allowing it to have good performance for noiseless recovery of
sparse signals, but in return displays a poor rate in the presence of
noise -- while min-$\ell_2$-norm interpolation has no structural bias
(except towards zero), causing it to fail to recover any non-zero signal even
in the absence of noise, but in return does not suffer from overfitting
of the noise. This behavior is also illustrated in
Figure~\ref{fig:tradeoff}.


%% file: sections/proof_sketch.tex
\ifarx
\section{Proof sketch}
\else
\section{PROOF SKETCH}
\fi

\label{sec:proof_sketch}

\label{sec:proofsketch}
In this section we present the main ingredients that are key to prove our risk upper bound of 
$\frac{\sigma^2}{\log(d/n)} + O \left( \frac{\sigma^2}{\log(d/n)^{3/2}} \right)$.
The proof sketch is interleaved with remarks providing insights on how our technique allows to improve upon previous works.
For the sake of clarity, we omit the discussion of the matching lower bound, as its proof follows exactly the same ideas.
The full proof is given in 
\ifarx Section\else Appendix\fi~\ref{sec:proof_maintext}.

The proof follows a standard localization/uniform convergence argument, where we first upper-bound the $\ell_1$-error  $ \|\what - \wstar\|_1 $ (localization) and then uniformly upper-bound the risk (i.e.\ $\ell_2$-error) over all interpolators with bounded $\ell_1$-error (uniform convergence).





\subsection{Localization}
\label{subsec:l1proofsketch}
We derive a high-probability upper bound on the $\ell_1$-error: $\norm{\hw-\wstar}_1 \leq B(n,d)$,
implying that the estimator of interest $\hw$ is an interpolator located in the $\ell_1$-ball of radius $B(n,d)$ centered at $\wstar$.

To upper-bound $\norm{\hw-\wstar}_1$ we first observe that, by definition of the BP estimator $\hw$ and via a simple triangle inequality \cite{chinot_2021}, it holds that


\begin{equation}
    \|\what - \wstar\|_1 \leq 2\sqrt{\sstar} \norm{\what - \wgt}_2 +  \min_{X w = \xi} \|w\|_1
\end{equation}
which is proven in
Lemma~\ref{lm:triangle_ineq}.
To control the first term, we make use of the loose high-probability upper bound $\norm{\what-\wgt}_2 \lesssim \sigma$ from previous works \cite{chinot_2021, wojtaszczyk_2010}. 
Thus we have $2\sqrt{\sstar} \norm{\what-\wgt}_2 \lesssim \sigma \sqrt{\sstar}$.
The second term, 
$\Phi_N := \min_{Xw = \xi} \|w\|_1$,
reflects how enforcing the interpolation of noise affects the $\ell_1$-norm of the estimator.
To control it with high probability directly is challenging, due to the randomness of both data $X$ and noise $\xi$.
Instead, we bound it using the Convex Gaussian Minimax Theorem (CGMT) \cite{thrampoulidis_2015}; we postpone the sketch of this derivation to Section~\ref{subsec:proofsketch_phin} where we show how we derive a high probability bound $\Phi_N \leq M(n,d)$ (the precise expression can be found in Proposition~\ref{prop:main_norm}). 
Having controlled the two terms separately, we get the high-probability bound
\ifarx
\begin{equation}
    \norm{\hw-w^*}_1
    \leq c \sigma \sqrt{\sstar} + M(n,d) =: B(n,d)
\end{equation}
\else
\begin{align}
    &\norm{\hw-w^*}_1 
    \leq \cgeoff \sigma \sqrt{\sstar} + M(n,d) =: B(n,d)
\end{align}
\fi
for some universal constant $\cgeoff > 0$.
Note that by assumption on $\sstar$, the first term is of order at most $\sqrt{\frac{\sigma^2 n}{\log(d/n)^5}}$, and hence negligible compared to $M(n,d) \geq \sqrt{\frac{\sigma^2 n}{2\log(d/n)}}$.

\paragraph{How tightness of $M(n,d)$ affects the upper bound and comparison to \cite{chinot_2021,wojtaszczyk_2010}.}
Intriguingly, our analysis requires a very precise expression for the deterministic upper bound on 
$\Phi_N = \min_{Xw=\xi} \norm{w}_1$.
Remember that $M(n,d)\approx \left(\frac{\sigma^2 n}{2 \log(d/n) - \log\log(d/n) - \log(\pi)}\right)^{1/2}$ is the precise upper bound used in our analysis, given in Proposition~\ref{prop:main_norm}.
We now discuss the repercussions of choosing a looser upper bound $\widetilde{M}(n,d)$ in our analysis:
 \begin{enumerate}
     \item for
     $\widetilde{M}(n,d) = M(n,d)\left(1 +  \frac{c}{\log(d/n)^2}\right)$ with $c>0$: We would still get exactly the precise upper bound of Theorem~\ref{thm:main}.
     \item for $\widetilde{M}(n,d) = M(n,d)\left(1 +  \frac{c}{\log(d/n)}\right)$ with $c>0$: We would obtain an upper bound of the correct order up to 
     a universal constant factor 
     $\norm{\what-\wgt}_2^2 \lesssim \frac{\sigma^2}{\log(d/n)}$.
     \item for $\widetilde{M}(n,d) = M(n,d)\left(1 + c\right)$ with $c>0$:
      The upper bound for $\Phi_N$ used in the analysis of  \cite{chinot_2021} is of order $\sqrt{\frac{2\sigma^2 n}{\log(d/n)}} \geq \widetilde{M}(n,d)$. Using this loose upper bound, our analysis would only yield an upper bound of constant order 
     $\norm{\what-\wgt}_2^2 \lesssim \sigma^2$, 
     same as obtained by \cite{chinot_2021,wojtaszczyk_2010}.
 \end{enumerate}



Yet, it is not clear whether the tightness of $M(n,d)$ is only needed due to our analysis, or whether any uniform convergence bound with the corresponding $B(n,d)$ would fail to yield tight bounds.
We leave this question as an interesting direction for future work, described further in Section~\ref{subsec:future_work}.



Note that the precise expression $M(n,d)$ alone
combined with previous analysis would only result in a constant upper bound of order $\sigma^2$ if inserted into the analysis 
conveyed in \cite{chinot_2021,wojtaszczyk_2010,koehler_2021} (discussed in the next section). 
In the following section we discuss the key steps that allow us to obtain the tight upper bound using our precise expression $M(n,d)$\footnote{Finally, we note that \cite{chinot_2021,wojtaszczyk_2010} both study the case where the noise can also be adversarial, for which the rate of order $\sigma^2$ for the risk is optimal.}. 

\subsection{Uniform risk bound and reduction to auxiliary problem by GMT}
Given that $\hw$ belongs with high probability to the set of interpolators located in the $\ell_1$-ball of radius $B(n,d)$ centered at $\wstar$, we proceed to upper-bound the risk of all such interpolators.

Concretely, we find a high-probability upper bound on
\ifarx
\begin{align}
    \Phi_+ :=& \max_{w} \norm{w-\wstar}_2^2
    \subjto \norm{w-w^*}_1 \leq B(n,d) ~\text{and}~ X(w-w^*) = \xi \\
    =& \max_w \norm{w}_2^2 \subjto \norm{w}_1 \leq B(n,d) ~\text{and}~ Xw = \xi.
\end{align}
\else
\begin{align}
    \Phi_+ :=& \max_{w} \norm{w-\wstar}_2^2 
    \subjto \norm{w-w^*}_1 \leq B(n,d)\\ 
    &\quad\quad\quad\quad\quad\quad\quad\quad ~\text{and}~ X(w-w^*) = \xi \\
    =& \max_w \norm{w}_2^2 \subjto \norm{w}_1 \leq B(n,d) ~\text{and}~ Xw = \xi.
\end{align}
\fi
While directly bounding $\Phi_+$ is challenging 
due to the randomness of both $X$ and $\xi$, we can instead make use of the Gaussian Minimax Theorem (GMT) which allows us to equivalently upper-bound the value of the so-called auxiliary problem
\ifarx
\begin{align}
\label{eq:proofsketch:phi+}
    \phi_+ := \max_w \norm{w}_2^2 \subjto \norm{w}_1 \leq B(n,d) ~\text{and}~ \innerprod{w}{h}^2 \geq (1-\rho) n (\sigma^2 + \norm{w}_2^2)
\end{align}
\else
\begin{align}
\label{eq:proofsketch:phi+}
    &\phi_+ := \max_w \norm{w}_2^2  \subjto \norm{w}_1 \leq B(n,d)\\
    &\quad\quad\quad\quad\quad ~\text{and}~ \innerprod{w}{h}^2 \geq (1-\rho) n (\sigma^2 + \norm{w}_2^2)
\end{align}
\fi
with $h$ an \iid\ Gaussian random vector
and $\cnoise$ a vanishing parameter.
Indeed, the GMT ensures that $\PP_{X, \xi}(\Phi_+ > t) \leq 2 \PP_h(\phi_+ \geq t) + \epsilon_{\cnoise}$ (see Proposition~\ref{prop:CGMT_application} for the expression of $\epsilon_{\cnoise}$);
in words, a high-probability upper bound on $\phi_+$ gives a high-probability upper bound on $\Phi_+$. 

Reformulating the optimization problem
\eqref{eq:proofsketch:phi+} as a one-dimensional optimization problem (see Section~\ref{subsec:pathapproach})
allows us to obtain the following tight upper bound 
\begin{equation}
    \phi_+ \leq \frac{\sigma^2}{\log(d/n)}\left(1 + \frac{c}{\log(d/n)^{1/2}}\right).
\end{equation}
We now discuss how a simpler relaxation of the constraints, used in previous work, leads to loose bounds.




\paragraph{Comparison to \cite{koehler_2021}.}
The optimization problem defining $\phi_+$ in Equation~\eqref{eq:proofsketch:phi+} is a maximization problem of a convex function over non-convex constraints.
We now show how to obtain a first loose upper bound 
following
\cite{koehler_2021} and briefly discuss why this methodology fails to give tight bounds (see also 
the paragraph
``Application: Isotropic features'' in that paper).

Using Hölder's inequality $\innerprod{w}{h} \leq \norm{w}_1 \norm{h}_{\infty}$, we obtain a proper relaxation of the problem~\eqref{eq:proofsketch:phi+} if we replace the constraints by 
\begin{equation}
    \locupper^2 \|h\|_\infty^2 \geq (1-\rho) n (\sigma^2 + \|w\|_2^2).
    \label{eq:koehlereq}
\end{equation} 
This immediately implies the upper bound
\begin{equation}
    \phidp \leq \frac{\gennormboundp(n,d)^2
\|\Hgaussian\|_{\infty}^2}{n(1-\cnoise)} - \sigma^2.
\end{equation}
However this bound is loose, even when we plug in our tight localization bound for $\locupper^2 \approx  \frac{\sigma^2 n}{2\log(d/n) - \log\log (d/n)}$. 
Indeed, with this estimate and by Gaussian concentration results,
the above bound reads for $d\gg n$
\ifarx
\begin{align}
    \phidp \leq 
    \frac{\locupper^2
    \|\Hgaussian\|_{\infty}^2}{n(1-\cnoise)} - \sigma^2  
    &\approx \frac{\sigma^2~ 2\log(d)}{2\log(d/n)} - \sigma^2 
    = \sigma^2 \frac{\log(n)}{\log(d/n)}.
    \label{eq:sometempeq}
\end{align}
\else
\begin{align}
    \phidp &\leq 
    \frac{\locupper^2
    \|\Hgaussian\|_{\infty}^2}{n(1-\cnoise)} - \sigma^2  \\
    &\approx \frac{\sigma^2~ 2\log(d)}{2\log(d/n)} - \sigma^2 
    = \sigma^2 \frac{\log(n)}{\log(d/n)}.
        \label{eq:sometempeq}
\end{align}
\fi
Note that this bound is constant in any polynomial growth regime $d \asymp n^\beta$, while we prove an upper bound in Theorem~\ref{thm:main} which vanishes in these regimes as $d,n \to \infty$. 

This looseness points to the fact that using H\"olders inequality is too imprecise.
In the next subsection, we describe our refined analysis which better takes into account the relationship between $\innerprod{w}{\Hgaussian}, \|w\|_2$ and $\|w\|_1$.

\subsection{Path approach: reparametrizing the auxiliary problem as a one-dimensional problem}
\label{subsec:pathapproach}

The key observation that allows us to derive a tight bound for $\phi_+$,
is that we can cast the $d$-dimensional problem \eqref{eq:proofsketch:phi+}
into a one-dimensional problem,
which we can study explicitly.
Namely, we identify a path $\gamma: \RR \to \RR^d$ for which we show that the optimum 
in~\eqref{eq:proofsketch:phi+}
is necessarily attained at $B(n,d) \gamma(\alpha)/\norm{\gamma(\alpha)}_1$ for some $\alpha \in \RR$.
Note that our reduction is exact, 
not a relaxation.

More precisely, we define the path
$\gamma: [1, \alpham] \to \RR^d$ by
\begin{equation}
    \gamma(\alpha) = \argmin_w \norm{w}_2^2
    \subjto \begin{cases}
        \innerprod{w}{h} = \norm{h}_\infty \\
        \forall i, h_i w_i \geq 0 \\
        \norm{w}_1 = \alpha
    \end{cases}
\end{equation}
(in particular $\norm{\gamma(\alpha)}_1=\alpha$)
and we show that 
\ifarx
\begin{align}
\label{eq:proof_sketch:pathformulation}
    \phi_+ = \locupper^2 \max_{1\leq \alpha\leq \alpha_{\max}}  \left(\frac{\|\gamma(\alpha)\|_2}{\alpha}\right)^2 \:\: \subjto \:\: &\frac{\|H\|^2_\infty}{(1-\rho)n } \geq \frac{\alpha^2\sigma^2}{\locupper^2} + \|\gamma(\alpha)\|_2^2
\end{align}
\else
\begin{align}
\label{eq:proof_sketch:pathformulation}
    \phi_+ &= \locupper^2 \max_{1\leq \alpha\leq \alpha_{\max}}  \left(\frac{\|\gamma(\alpha)\|_2}{\alpha}\right)^2 \:\: 
    \\\text{ s.t. } \:\: &\frac{\|h \|^2_\infty}{(1-\rho)n} \geq \frac{\alpha^2\sigma^2}{\locupper^2} + \|\gamma(\alpha)\|_2^2
\end{align}
\fi
(see \ifarx Section\else Appendix\fi~\ref{subsec:proof_of_props_phi}).
Because $\gamma(\alpha)$ is the argmin of a convex optimization problem, it is relatively easy to study, and we can even derive an exact expression for it (Lemma~\ref{lm:charact_gamma}). We now discuss the two key steps to study the optimization problem in Equation~\eqref{eq:proof_sketch:pathformulation}.






\paragraph{a) Monotonicity of the objective.} We observe that $\frac{\|\gamma(\alpha)\|_2}{\alpha}$ is monotonically decreasing and that $\norm{\gamma(\alpha)}_2^2$ is a convex function (Lemma~\ref{lm:interval}). 
This has two important consequences.
Firstly, the set of $\alpha$'s which satisfy the constraints in~\eqref{eq:proof_sketch:pathformulation} is an interval, denoted $[\ulalpha_I, \olalpha_I]$. 
Secondly, denoting
$\alpha^*$ ($\in [\ulalpha_I, \olalpha_I]$)
an argmax of~\eqref{eq:proof_sketch:pathformulation},
we have for any $\alpha < \ulalpha_I$
\begin{equation}
\label{eq:upperboundsimple}
    \locupper^2 \left(\frac{\|\gamma(\alpha)\|_2}{\alpha}\right)^2 \geq \locupper^2 \left(\frac{\|\gamma(\alpha^*)\|_2}{\alpha^*}\right)^2 = \phi_+.
\end{equation}
So to obtain an upper bound on $\phi_+$, all we need is to find an $\alpha$ on the left of the feasible interval $[\ulalpha_I, \olalpha_I]$.
We do this by finding both an $\alpha_n$ such that ${ \ulalpha_I \leq \alpha_n \leq \olalpha_I }$ (i.e.\ $\alpha_n$ is feasible), and an $\alpha_s$ such that $\alpha_s < \alpha_n$ and $\alpha_s \not\in [\ulalpha_I, \olalpha_I]$ (i.e.\ $\alpha_s$ is not feasible).




\paragraph{b) Discretization of the path.}
We observe that there exist ``breakpoints'' $1 = \alpha_2 < ... < \alpha_{d+1} = \alphamax$ for which $\gamma(\alpha_s)$ has a special structure (in particular it is $(s-1)$-sparse).
Further, applying Gaussian concentration results to $h$ leads to high-probability estimates for 
$\alpha_s$
and $\norm{\gamma(\alpha_s)}_2$
(Proposition~\ref{prop:concentration_gammas}).
Thanks to those estimates, 
we show that $\alpha_n$ is feasible for \eqref{eq:proof_sketch:pathformulation}
and
we find a choice of $s < n$ such that $\alpha_s$ is not feasible, with high probability.
Thus, with high probability $\phi_+$ is upper-bounded by $\locupper^2 \left(\frac{\|\gamma(\alpha_s)\|}{\alpha_s}\right)^2$ -- for which we have high-probability estimates.

\paragraph{Intuition for the definition of $\gamma(\alpha)$.}
As discussed in the previous subsection, the relaxation of \eqref{eq:phidp} based on H\"older's inequality $\innerprod{w}{h} \leq \norm{w}_1 \norm{h}_\infty$, used
\ifarx in the paper \cite{koehler_2021},
\else by \citet{koehler_2021},
\fi 
is too loose.
Informally, it effectively amounts to forgetting the direction of $h$ and only optimizing over the $\ell_1$ and $\ell_2$-norms of vectors.
The main idea of our refined analysis is to introduce a path, $\{\gamma(\alpha)/\alpha\}_\alpha$, allowing us to better take into account the relationship between $\innerprod{w}{\Hgaussian}, \|w\|_2$ and $\|w\|_1$.
To intuitively understand how the path achieves this goal, it may be easier to use its following form:
\begin{equation} 
    \olw(\beta) = \argmax_w \innerprod{w}{h} \subjto \begin{cases}
        \norm{w}_2^2 \leq \beta \\
        \forall i, h_i w_i \geq 0 \\
        \norm{w}_1 = 1
    \end{cases},
\end{equation}
which is equivalent to $\gamma(\alpha)$ up to linear reparametrization and rescaling (see Equation~\eqref{eq:def_olw_beta} in Appendix~\ref{apx:subsec:param_argmaxmin}).
Note that $\langle \olw(1), h \rangle = \norm{h}_{\infty} \norm{\olw(1)}_1$, which exactly recovers the equality case of H\"older's inequality (i.e.\ $\olw(1)$ is a subgradient of the $\ell_1$-norm at $h$).
More generally, $\langle w,h\rangle \leq \langle \olw(\beta),h\rangle \norm{w}_1$ for any $w$ such that $\norm{w}_2^2 \leq \beta \norm{w}_1^2$, which can be understood as a refined H\"older's inequality for limited $\ell_2$-norms.

\subsection{Obtaining a good estimate for \texorpdfstring{$\Phi_N$}{PhiN}}
\label{subsec:proofsketch_phin}
Finally, we unveil how we derive the high-probability upper bound $\Phi_N \leq M(n,d)$ in the localization step (1.):
The derivation actually uses the same tools as for the upper bound of $\Phi_+$.
Using the Convex Gaussian Minimax Theorem (CGMT) \cite{thrampoulidis_2015},
which is a variant of the GMT for convex-concave functions,
we can again introduce an auxiliary problem
\begin{equation}
\label{eq:proofsketch:phin}
    \phi_N := \min_w \norm{w}_1 \subjto \innerprod{w}{h}^2 \geq (1+\rho) n (\sigma^2 + \|w\|_2^2)
\end{equation}
with $h$ an i.i.d.\ Gaussian vector and $\rho$ a vanishing parameter,
with the property that high-probability upper bounds on $\phi_N$ give high-probability upper bounds on $\Phi_N$.
Further,
we can again reduce this $d$-dimensional optimization problem to one over the same path $\{\gamma(\alpha)\}_\alpha$:
\begin{equation}
\label{eq:constarintsphin}
    \phi_N = \min_{1\leq \alpha\leq\alphamax}  f(\alpha)
    \subjto
    \norm{h}_\infty^2 \geq (1+\cnoise) n \norm{\gamma(\alpha)}_2^2
\end{equation}
(see \ifarx Section\else Appendix\fi~\ref{subsec:CGMTproof} for the expression of $f(\alpha)$).

Since we want to upper-bound this minimum, it is sufficient to find some $\alpha$ which satisfies the constraints in Equation~\eqref{eq:constarintsphin}.
In particular, we again 
focus on
the breakpoints $\{\alpha_s\}_{s \in \{2,...,d+1\}}$,
and show that with high probability
$\alpha = \alpha_n$ is
a valid choice which is approximately tight (see Remark~\ref{rk:choice_s_n}).

\paragraph{Sparsity of  $\gamma(\alpha_n)$.}
In summary, the proof is essentially based on the localization around a rescaled version of $\gamma(\alpha_n)$, which is a $(n-1)$-sparse vector. 
This choice can also be motivated by a different argument:
It is well known that the minimizer of the optimization problem
$\min_{Xw=\xi} \norm{w}_1$
defining  $\PhiN$ is $n$-sparse. Hence, due to the strong 
connection between
the optimization problems defining $\phidn$ and $\PhiN$, we also expect the minimizer of Equation~\eqref{eq:constarintsphin} to be approximately $n$-sparse. 

%% file: sections/proof_maintext.tex
\ifarx
\section{Proof of main result} 
\else
\section{PROOF OF MAIN RESULT} 
\fi

\label{sec:proof_maintext}
In this section, we present the proof of our main result, Theorem~\ref{thm:main}. 
In \ifarx Section~\ref{subsec:proof_of_thm}\else Appendix~\ref{subsec:proof_of_thm}\fi, we describe the main steps of the proof rigorously, in the form of three propositions which we then prove in \ifarx Section\else Appendix\fi~\ref{subsec:proof_of_props_phi}.
Full proofs for the intermediary Lemmas and Propositions are given in Appendix~\ref{apx:sec:appendix_proofs}.

\paragraph{Notation.}
On the finite-dimensional space $\RR^d$, we write $\norm{\cdot}_2$ for the Euclidean norm and $\innerprod{\cdot}{\cdot}$ for the Euclidean inner product.
The $\ell_1$ and $\ell_{\infty}$-norms are denoted by $\norm{\cdot}_1$ and $\norm{\cdot}_\infty$, respectively.
The vectors of the standard basis are denoted by $e_1, ..., e_d$, and
$\onevec \in \RR^d$ is the vector with all components equal to $1$.
For $s \leq d$ and $H \in \RR^d$, $\Hvecs$ is the vector such that $(\Hvecs)_i = H_i$ if $i \leq s$ and $0$ otherwise.
$\NNN(\mu, \Sigma)$ is the normal distribution with mean $\mu$ and covariance $\Sigma$, $\Phi(x)$ is
the cumulative distribution function of the scalar standard normal distribution, $\Phic(x) = 1-\Phi(x)$,
and $\log$ denotes the natural logarithm.
For all $s \leq d$, we denote by $t_s \in \RR$ the quantile of the standard normal distribution defined by
$2\Phic(t_s) = s/d$.
The $n$ samples $x_i \in \RR^d$
form the rows of the data matrix \mbox{$X = \left[ x_1 ~...~ x_n
    \right]^\top$}, with $X_{ij} \sim \NNN(0,1)$ for each $i,j$.  The
scalars $y_i$, $\xi_i$ are also aggregated into vectors $y, \xi \in
\RR^n$ with $\xi \sim \NNN(0, \sigma^2 I_n)$ and $y = Xw^* + \xi$.
With this notation, $\hw$ interpolates the data $X\hw = y$ which is equivalent to
$X(\hw-w^*) = \xi$.
To easily keep track of the dependency on dimension and
sample size,
we reserve the $O(\cdot)$ notation to contain only universal constants, without any hidden dependency
on $d$, $n$, $\sigma^2$ or $\norm{w^*}_0$.
We will also use $c_1,c_2,...$ and $\kappa_1,\kappa_2,...$ to denote positive universal constants reintroduced each time in the proposition and lemma statements, except for $\cgeoff$ and $\cbrho$ which should be considered as fixed throughout the whole proof.

\subsection{Proof of Theorem~\ref{thm:main}}
\label{subsec:proof_of_thm}

We proceed by a localized uniform convergence approach,
similar to 
\ifarx the papers
\else \fi
\cite{chinot_2021,koehler_2021,ju_2020,muthukumar_2020}, and common in the literature, \eg, on structural risk minimization.
That is, the proof consists of two steps:
\begin{enumerate}
    \item \emph{Localization.}
    We prove that, with high probability, the min-$\ell_1$-norm interpolator $\hw$ satisfies
    $\norm{\hw-w^*}_1 \leq c \sigma \sqrt{\sstar} + \min_{Xw=\xi} \norm{w}_1$ for some universal constant $c>0$. 
    We then derive a (finer than previously known) high-probability upper bound on the second term,
    \begin{equation} \label{eq:Phi_N} \tag{$P_N$}
        \min_{Xw=\xi} \norm{w}_1 
        ~ =: \Phi_N 
        ~ \leq M(n,d).
    \end{equation}
    Consequently, with high probability $\hw$ satisfies
    \begin{equation}
        \norm{\hw-w^*}_1 \leq c \sigma \sqrt{\sstar} + M(n,d) =: B(n,d).
    \end{equation}
    
    \item \emph{Uniform convergence.}
    We derive high-probability uniform upper and lower bounds on the prediction error for all interpolators located no farther than $B(n,d)$ from $w^*$ in $\ell_1$ norm.
    In symbols,
    we find a high-probability upper bound for
\ifarxiv
    \begin{align} \label{eq:Phi_+} \tag{$P_+$}
        \max_{\substack{
            \norm{w-w^*}_1 \leq B(n,d) \\
             X(w-w^*)=\xi
        }} \norm{w-w^*}_2^2 
        = \max_{\substack{
            \norm{w}_1 \leq B(n,d) \\
             Xw=\xi
        }} \norm{w}_2^2 
        ~
        =: \Phi_+
    \end{align}
\else
    \begin{align} \label{eq:Phi_+} \tag{$P_+$}
        \max_{\substack{
            \norm{w-w^*}_1 \leq B(n,d) \\
             X(w-w^*)=\xi
        }} & \norm{w-w^*}_2^2 \\
        = \max_{\substack{
            \norm{w}_1 \leq B(n,d) \\
             Xw=\xi
        }} & \norm{w}_2^2 
        ~
        =: \Phi_+
    \end{align}
\fi
    and a high-probability lower bound for
\ifarxiv
    \begin{align} \label{eq:Phi_-} \tag{$P_-$}
        \min_{\substack{
            \norm{w-w^*}_1 \leq B(n,d) \\
             X(w-w^*)=\xi
        }} \norm{w-w^*}_2^2 
        = \min_{\substack{
            \norm{w}_1 \leq B(n,d) \\
             Xw=\xi
        }} \norm{w}_2^2 
        ~
        =: \Phi_-.
    \end{align}
\else
    \begin{align} \label{eq:Phi_-} \tag{$P_-$}
        \min_{\substack{
            \norm{w-w^*}_1 \leq B(n,d) \\
             X(w-w^*)=\xi
        }} & \norm{w-w^*}_2^2  \\
        = \min_{\substack{
            \norm{w}_1 \leq B(n,d) \\
             Xw=\xi
        }} & \norm{w}_2^2 
        ~
        =: \Phi_-.
    \end{align}
\fi
\end{enumerate}
By definition of $B(n,d)$ in \eqref{eq:Phi_N}, with high probability the min-$\ell_1$-norm interpolator $\what$ belongs to the set of feasible solutions in \eqref{eq:Phi_+} and \eqref{eq:Phi_-},
and hence the second step yields high-probability upper and lower bounds on its prediction error $\norm{\what-w^*}_2^2$. 

The key is thus to derive tight high-probability bounds for the quantities $\Phi_N, \Phi_+, \Phi_-$.
Our derivation proceeds in two parts, described below. 
The first part uses the CGMT to convert the original optimization problem to an auxiliary problem, similar to \cite{koehler_2021}. 
The second part, which contains the crucial elements for our proof of the vanishing upper bound and is the key technical contribution of this paper, consists in reducing the $d$-dimensional auxiliary problem to a scalar one using a path reparametrization.

\paragraph{Preliminary: Localization around $w^*$.}

The following fact shows how, as announced, $\Phi_N$ can be used to derive a localization bound for $\hw$.
\begin{lemma} \label{lm:triangle_ineq}
    Suppose 
    ${\norm{\wgt}_0 \leq \cwone \frac{n}{\log(d/n)^{5}}}$ 
    for some universal constant $\cwone>0$.
    There exist universal constants $\cn,\cl,\cd, \cgeoff, \cprobgeoff>0$ such that, if
    $n \geq \cn$ and
    $\cl n \log(n)^2 \leq d \leq \exp(\cd n^{1/5})$,
    then
    the min-$\ell_1$-norm interpolator $\hw$ satisfies
    \begin{equation}
        \norm{\hw-w^*}_1 \leq \cgeoff \sigma \sqrt{\sstar} + \min_{Xw=\xi} \norm{w}_1
    \end{equation}
    with probability at least
    $1-d \exp\left( -\cprobgeoff n \right)$.
\end{lemma}
The proof of Lemma~\ref{lm:triangle_ineq} is given in Appendix~\ref{apx:subsec:triangle}.
Interestingly, it makes use of the loose upper bound $\norm{\hw-w^*}_2^2 \lesssim \sigma^2$, shown previously by \cite{wojtaszczyk_2010} and \cite{chinot_2021}, as an intermediate result.

\subsubsection{(Convex) Gaussian Minimax Theorem}
\label{subsec:CGMTproof}
Since each of the quantities $\Phi_N, \Phi_+, \Phi_-$ is defined as the optimal value of a stochastic program with Gaussian parameters, we may apply the (Convex) Gaussian Minimax Theorem ((C)GMT)
\cite{gordon_1988,thrampoulidis_2015}.
On a high level, given a ``primary'' optimization program with Gaussian parameters, the (C)GMT relates it to an ``auxiliary'' optimization program,
so that high-probability bounds on the latter imply high-probability bounds on the former. The following proposition applies the CGMT on $\Phi_N$ and the GMT on $\Phi_+$, $\Phi_-$.

\begin{proposition} 
\label{prop:CGMT_application}
    For $\Hgaussian \sim \NNN(0,I_d)$, define the stochastic auxiliary optimization problems:
\ifarxiv
    \begin{align}
    \tag{$A_N$} \label{eq:phidn}
        \phidn(\cnoise) &= \min_w \norm{w}_1  
        \subjto \innerprod{w}{\Hgaussian}^2 \geq (1+\cnoise) n (\sigma^2 + \norm{w}_2^2) \\
    \tag{$A_+$} \label{eq:phidp}
        \phidp(\cnoise) &= \max_w \norm{w}_2^2 
        \subjto \begin{cases}
            \norm{w}_1 \leq B(n,d) \\
            \innerprod{w}{\Hgaussian}^2 \geq (1-\cnoise) n (\sigma^2 + \norm{w}_2^2)
        \end{cases} \\
    \tag{$A_-$} \label{eq:phidm}
        \phidm(\cnoise) &= \min_w \norm{w}_2^2 
        \subjto  \begin{cases}
            \norm{w}_1 \leq B(n,d) \\
            \innerprod{w}{\Hgaussian}^2 \geq (1-\cnoise) n (\sigma^2 + \norm{w}_2^2)
        \end{cases}
    \end{align}
\else
    \begin{align}
    \tag{$A_N$} \label{eq:phidn}
        \phidn(\cnoise) &= \min_w \norm{w}_1  \\
        &\subjto \innerprod{w}{\Hgaussian}^2 \geq (1+\cnoise) n (\sigma^2 + \norm{w}_2^2) \\
    \tag{$A_+$} \label{eq:phidp}
        \phidp(\cnoise) &= \max_w \norm{w}_2^2 \\
        &\subjto \begin{cases}
            \norm{w}_1 \leq B(n,d) \\
            \innerprod{w}{\Hgaussian}^2 \geq (1-\cnoise) n (\sigma^2 + \norm{w}_2^2)
        \end{cases} \\
    \tag{$A_-$} \label{eq:phidm}
        \phidm(\cnoise) &= \min_w \norm{w}_2^2 \\
        &\subjto  \begin{cases}
            \norm{w}_1 \leq B(n,d) \\
            \innerprod{w}{\Hgaussian}^2 \geq (1-\cnoise) n (\sigma^2 + \norm{w}_2^2)
        \end{cases}
    \end{align}
\fi
    where $0<\cnoise<1/2$ can be any small enough quantity.
    For any $t \in \RR$, it holds that
    \begin{align}
        \PP( \Phi_N > t ) &\leq 2\PP( \phidn(\cnoise) \geq t ) + 6 \exp\left(-\frac{n \cnoise^2}{100} \right) \label{eq:probeq1} \\
        \text{and}~~~ \PP( \Phi_+ > t ) &\leq 2\PP( \phidp(\cnoise) \geq t ) + 6 \exp\left(-\frac{n \cnoise^2}{100}\right)   \label{eq:probeq2} \\
        \text{and}~~~   \PP( \Phi_- < t ) &\leq 2\PP( \phidm(\cnoise) \leq t ) + 6 \exp\left(-\frac{n \cnoise^2}{100}\right) \label{eq:probeq3} ,
    \end{align}
    where on the left-hand side $\PP$ denotes the probability distribution  over $X$ and $\xi$, and on the right-hand side the distribution over $\Hgaussian$.
\end{proposition}

For the remainder of this proof, we choose%
\footnote{
    This choice of $\cnoise$ is justified by the proof of Proposition~\ref{prop:main_norm}. Indeed, for an arbitrary choice of $\cnoise < 1/2$, one could still show the same bound with just an extra factor:
    $(\phidn)^2 \leq (1+\cnoise) \Brho(n,d)$,
    holding with still the same probability.
    This would translate to a bound on $\Phi_N$ holding with probability 
    ${1 - 12 \exp\left(-2\frac{n}{\log(d/n)^5}\right)} - 6 \exp\left(-\frac{n \cnoise^2}{100} \right)$.
    So the choice $\cnoise = \frac{10}{\log(d/n)^{5/2}}$ ``comes at no cost'' in terms of the probability with which the bound holds,
    while being sufficiently small to allow for a satisfactory bound (it only affects the constant $\cbrho$ appearing in $\Brho(n,d)$).
}
\begin{equation}
    \cnoise = \frac{10}{\log(d/n)^{5/2}}.
\end{equation}
As such, from now on, we simply write $\phi_N$, $\phi_+$, $\phi_-$. 
The proof of Proposition~\ref{prop:CGMT_application}, given in Appendix~\ref{apx:subsec:CGMT}, closely follows Lemmas~{3-7} in the paper \cite{koehler_2021}.
For clarity, note that the three pairs of stochastic programs 
(\ref{eq:Phi_N}/\ref{eq:phidn}),
(\ref{eq:Phi_+}/\ref{eq:phidp}),
(\ref{eq:Phi_-}/\ref{eq:phidm})
are not coupled: Proposition~\ref{prop:CGMT_application} should be understood as consisting of three separate statements, each using a different independent copy of $h$.

As a result of the proposition, the goal of finding high-probability bounds on $\Phi_N, \Phi_+, \Phi_-$ now reduces to finding high-probability bounds on $\phi_N$, $\phi_+$, $\phi_-$, respectively.

\subsubsection{Bounds on \texorpdfstring{$\phidn, \phidp, \phidm$}{phiN,phi+,phi-}}

To obtain tight bounds on the auxiliary quantities $\phidn, \phidp, \phidm$, we adopt a significantly different approach from previous works.
The main idea is to reduce the
optimization problems \eqref{eq:phidn}, \eqref{eq:phidp} and \eqref{eq:phidm} to optimization problems over a parametric 
path~$\{\gamma(\alpha)\}_{\alpha} \subset \RR^d$.
Here we only state the results and refer to \ifarx Section\else Appendix\fi~\ref{subsec:proof_of_props_phi} for their proofs and further intuition.
For the remainder of this proof,
we denote by $t_n \in \RR$ the quantile of the standard normal distribution defined by
$2\Phic(t_n) = n/d$.

\begin{proposition} 
\label{prop:main_norm}
    There exist universal constants $\cn,\cl,\cd, \cbrho>0$ such that, if
    $n \geq \cn$ and
    $\cl n \leq d \leq \exp(\cd n^{1/5})$,
    then
    \begin{align}
        \phidn^2 &\leq 
        \frac{\sigma^2 n}{t_n^2} \left(
            1 - \frac{2}{t_n^2} + \frac{\cbrho}{t_n^4}
        \right)
    \end{align}
    with probability at least ${1 - 6 \exp\left(-2\frac{n}{\log(d/n)^5}\right)}$ over the draws of $\Hgaussian$. 
    
    Consequently, \eqref{eq:Phi_N} holds with
    \begin{equation} \label{eq:def_Mnd}
        M(n,d) := \sqrt{
            \frac{\sigma^2 n}{t_n^2} \left(
                1 - \frac{2}{t_n^2} + \frac{\cbrho}{t_n^4}
            \right)
        }
    \end{equation}
    with probability at least $1-18 \exp\left( -\frac{n}{\log(d/n)^5} \right)$
    over the draws of $X$ and $\xi$.
\end{proposition}

Hence by Lemma~\ref{lm:triangle_ineq}, the min-$\ell_1$-norm interpolator is located close to the true vector $w^*$, namely the $\ell_1$ distance is bounded by the deterministic quantity
\begin{align}
\label{eq:bnequaiton}
    \norm{\hw-w^*}_1 \leq c \sigma \sqrt{\sstar} + M(n,d) =: B(n,d)
\end{align}
 with probability at least
$1 - 18 \exp\left(-\frac{n}{\log(d/n)^5}\right) - d \exp\left( -\cprobgeoff n \right)$ and with $\cgeoff, \cprobgeoff>0$ some universal constants. 
We now establish high-probability upper resp.\ lower bounds for $\phidp$ resp.\ $\phidm$.

\begin{proposition}
\label{prop:main_upper_lower}
    Suppose 
    ${\norm{\wgt}_0 \leq \cwone \frac{n}{\log(d/n)^{5}}}$ 
    for some universal constant $\cwone>0$.
    There exist universal constants $\cn,\cl,\cd,\cbound, \cprobgeoff>0$ such that, if
    $n \geq \cn$ and
    $\cl n \leq d \leq \exp(\cd n^{1/5})$,
    then 
    each of the two events
\ifarxiv
    \begin{align}
        \phidp &\leq \frac{\sigma^2}{\log(d/n)} \left(1+ \frac{\cbound}{\sqrt{\log(d/n)}}\right)
        & &\text{and} &
        \phidm &\geq \frac{\sigma^2}{\log(d/n)} \left(1 - \frac{\cbound}{\sqrt{\log(d/n)}}\right)
    \end{align}
\else
    \begin{align}
        \phidm &\geq \frac{\sigma^2}{2\log(d/n)} \left(1- \frac{\cbound}{\sqrt{\log(d/n)}}\right) \\
        \text{and}~
        \phidp &\leq \frac{\sigma^2}{\log(d/n)} \left(1+ \frac{\cbound}{\sqrt{\log(d/n)}}\right)
    \end{align}
\fi
    happens with probability at least
    ${ 1 - 18 \exp\left(-\frac{n}{\log(d/n)^5}\right) }$
    over the draws of $\Hgaussian$.
\end{proposition}

Theorem~\ref{thm:main} follows straightforwardly from
Lemma~\ref{lm:triangle_ineq} and Propositions~\ref{prop:CGMT_application}, \ref{prop:main_norm} and \ref{prop:main_upper_lower}.

\subsection{Proof of Propositions~\ref{prop:main_norm} and \ref{prop:main_upper_lower}}
\label{subsec:proof_of_props_phi}

In this section we detail our analysis of the auxiliary optimization problems \eqref{eq:phidn}, \eqref{eq:phidp} and \eqref{eq:phidm}. 
We start by a remark that considerably simplifies notation:
The definitions of $\phidn, \phidp, \phidm$ are unchanged if,
in \eqref{eq:phidn}, \eqref{eq:phidp}, \eqref{eq:phidm},
$\Hgaussian$ is replaced by the reordered vector of its absolute order statistics, \ie, by $H$ such that $H_i$ is the $i$-th largest absolute value of $h$.
Throughout this proof, we condition on the event where $H$ has distinct and positive components:
\mbox{$H_1 > ... > H_d > 0$}, which holds with probability one. 
Henceforth, unless specified otherwise, references to the optimization problems~\eqref{eq:phidn}, \eqref{eq:phidp} and \eqref{eq:phidm} refer to the equivalent problems where $\Hgaussian$ is replaced by $H$.
Also recall that
we choose
$\cnoise = \frac{10}{\log(d/n)^{5/2}}$.
The key steps in the proof of Propositions~\ref{prop:main_norm} and \ref{prop:main_upper_lower} are as follows.
\begin{itemize}
    \item For each of the three optimization problems \eqref{eq:phidn}, \eqref{eq:phidp} and \eqref{eq:phidm}, 
    we show that the argmax (or argmin) is of the form $\av \gamma(\alpha)$ for some $\av>0$ and a parametric path
    $\Gamma = \lbrace\frac{\gamma(\alpha)}{\alpha}\rbrace_\alpha$
    (which depends on $H$).
    Hence we can restate \eqref{eq:phidn}, \eqref{eq:phidp} and \eqref{eq:phidm} as optimization problems over
    a scalar variable $\alpha$ and a scale variable $\av>0$. 
    (\ifarx Section\else Appendix\fi~\ref{subsubsec:parametrizing_argmaxmin})
    
    \item Still conditioning on $H$, we explicitly characterize the parametric path $\Gamma$. 
    In particular, we show that it is piecewise linear with breakpoints $\gamma(\alphathr{s})$ having closed-form expressions.
    (\ifarx Section\else Appendix\fi~\ref{subsubsec:characterizing_path})
    
    \item Thanks to the concentration properties of $H$ (\ifarx Section\else Appendix\fi~\ref{subsubsec:concentration_gammas}),
    evaluating at one of the breakpoints yields the desired high-probability upper bound on $\phidn$ (\ifarx Section\else Appendix\fi~\ref{subsubsec:proof_bound_phidn}).
    
    \item A fine-grained study of the intersection of 
    $\RR_+ \Gamma := \lbrace \av\frac{\gamma(\alpha)}{\alpha} \rbrace_{b \in \RR_+, \alpha}$
    with the constraint set of \eqref{eq:phidp} and \eqref{eq:phidm}, as well as the concentration properties of $H$, yield the desired high-probability bounds on $\phidp$ and $\phidm$.
    (\ifarx Section\else Appendix\fi~\ref{subsubsec:proof_bound_phidp_phidm})
\end{itemize}


\subsubsection{Parametrizing the argmax/argmin}
\label{subsubsec:parametrizing_argmaxmin}

Note that in the optimization problems \eqref{eq:phidn}, \eqref{eq:phidp} and \eqref{eq:phidm},
the variable $w$ only appears through $\norm{w}_2$, $\norm{w}_1$ and $\innerprod{w}{H}$. Thus, we can add the constraint that $\forall i,~w_i \geq 0$ without affecting the optimal solution.  
We will show that the path $\Gamma = \left\lbrace\frac{\gamma(\alpha)}{\alpha}\right\rbrace_\alpha$ can be used to parametrize the solutions of the optimization problems, where
$\gamma:  [1, \alpham] \to \RR^d$ is defined by
\begin{equation}
    \gamma(\alpha) = \argmin_w \norm{w}_2^2
    \subjto \begin{cases}
        \innerprod{w}{H} \geq \norm{H}_\infty \\
        \forall i, w_i \geq 0 \\
        \onevec^\top w = \norm{w}_1 = \alpha
    \end{cases}
\end{equation}
and $\alpham = d \frac{\norm{H}_\infty}{\norm{H}_1}$.
Specifically, the following key lemma states that (at least one element of) the argmax/argmin of \eqref{eq:phidn}, \eqref{eq:phidp} and \eqref{eq:phidm}
is of the form $\av \frac{\gamma(\alpha)}{\alpha}$ for some $\av > 0$ and $\alpha \in [1, \alpham]$.
This allows to reduce the optimization problems to a single scalar variable and a scale variable. 

\begin{lemma}
\label{lm:parametrization}
    Denoting for concision $\gennormboundp = \gennormboundp(n,d)$,
    we have that:
\begin{enumerate}
    \item The variable $w$ in \eqref{eq:phidn} can equivalently be constrained to belong to the set $\RR_+ \Gamma$, i.e., 
    \begin{align}
    \label{eq:phidn_reparam} \tag{$A_N'$}
        \phidn 
        =
        \min_{\av>0, 1 \leq \alpha \leq \alpham} \av 
        \subjto~
        \av^2 \norm{H}_\infty^2 \geq (1+\cnoise) n (\sigma^2 \norm{\gamma(\alpha)}_1^2 + \av^2 \norm{\gamma(\alpha)}_2^2).
    \end{align}
    
    \item The variable $w$ in \eqref{eq:phidp} can equivalently be constrained to belong to the set $B\Gamma$, i.e., 
    \begin{align}
    \label{eq:phidp_reparam} \tag{$A_+'$}
        \phidp
        =
        \max_{1 \leq \alpha \leq \alpham} \gennormboundp^2 \frac{\norm{\gamma(\alpha)}_2^2}{\norm{\gamma(\alpha)}_1^2} 
        \subjto~
        \gennormboundp^2 \norm{H}_\infty^2 \geq (1-\cnoise) n (\sigma^2 \norm{\gamma(\alpha)}_1^2 + \gennormboundp^2 \norm{\gamma(\alpha)}_2^2).
    \end{align}
    
    \item The variable $w$ in \eqref{eq:phidm} can equivalently be constrained to belong to the set $(0,B] \Gamma$, i.e., 
    \begin{align}
    \label{eq:phidm_reparam} \tag{$A_-'$}
        &\phidm
        =
        \min_{\substack{
            0<\av \leq \gennormboundp \\
            1 \leq \alpha \leq \alpham
        }}
        \av^2 \frac{\norm{\gamma(\alpha)}_2^2}{\norm{\gamma(\alpha)}_1^2} 
        \subjto~ 
        \av^2 \norm{H}_\infty^2 \geq (1-\cnoise) n (\sigma^2 \norm{\gamma(\alpha)}_1^2 + \av^2 \norm{\gamma(\alpha)}_2^2).
    \end{align}
\end{enumerate}
\end{lemma}

The proof of the lemma is given in Appendix~\ref{apx:subsec:param_argmaxmin}.
To give an intuitive explanation for the equivalence between \eqref{eq:phidn} and \eqref{eq:phidn_reparam}, consider a penalized version of \eqref{eq:phidn}:
$\min_w \norm{w}_1 - \lambda \left(
    \innerprod{w}{\Hgaussian}^2
    -
    (1+\cnoise) n (\sigma^2 + \norm{w}_2^2)
\right)$
with $\lambda>0$.
For fixed values of $\norm{w}_1$ and $\innerprod{w}{\Hgaussian}$, minimizing this penalized objective is equivalent to minimizing $\norm{w}_2^2$.
Hence, we can expect the argmin to be attained at $\av \frac{\gamma(\alpha)}{\alpha}$ for some $\av>0,\alpha$. 

\subsubsection{Characterizing the parametric path}
\label{subsubsec:characterizing_path}

As $\gamma(\alpha)$ is defined as the optimal solution of a convex optimization problem,
we are able to obtain a closed-form expression, by a straightforward application of Lagrangian duality.
The only other non-trivial ingredient is to notice that, at optimality, the inequality constraint $\innerprod{w}{H} \geq \norm{H}_\infty$
necessarily holds with equality.
Denote $\Hvecs$ the vector equal to $H$ on the first $s$ components and $0$ on the last $(d-s)$, and similarly for $\onevecs$.
Define, for any integer $2 \leq s \leq d$,
\begin{align}
     \alphathr{s} =
    \frac{
        \left( \norm{\Hvecs}_1 - s H_s \right) \norm{H}_\infty
    }{
        \norm{\Hvecs}_2^2 - \norm{\Hvecs}_1 H_s
    }.
\end{align}
Note that $\alphathr{2} = 1$.
Let $\alphathr{d+1} = \alpham$.
\begin{lemma}
\label{lm:charact_gamma}

    For all $1 < \alpha \leq \alpham$,
    denote $s$ the unique integer in $\{2,...,d\}$ such that
    \mbox{$\alphathr{s} < \alpha \leq \alphathr{s+1}$}.
    Then $\gamma(\alpha) = \lambda \Hvecs - \mu \onevecs$
    (in particular it is $s$-sparse)
    where the dual variables $\lambda$ and $\mu$ are given by
    \begin{align}
            \lambda &= \frac{1}{s \norm{\Hvecs}_2^2 - \norm{\Hvecs}_1^2}
            \left( s \norm{H}_\infty - \alpha \norm{\Hvecs}_1 \right) \\
            ~~~\text{and}~~~
            \mu &= \frac{1}{s \norm{\Hvecs}_2^2 - \norm{\Hvecs}_1^2}
            \left( \norm{\Hvecs}_1 \norm{H}_\infty - \alpha \norm{\Hvecs}_2^2 \right).
    \end{align}
\end{lemma}
The proof of the lemma is given in Appendix~\ref{apx:subsec:geometric_lemma}.

\subsubsection{Concentration of norms of  \texorpdfstring{$\gamma(\alphathr{s})$}{gamma(alphas))}}
\label{subsubsec:concentration_gammas}

Given the explicit characterization of the parametric path, we now
study its breakpoints $\gamma(\alphathr{s})$
(\mbox{$s \in \{2,...,d\}$}), 
and more precisely we estimate $\norm{\gamma(\alphathr{s})}_1 = \alphathr{s}$ and $\norm{\gamma(\alphathr{s})}_2$ as a function of $s$
(we have by definition $\innerprod{\gamma(\alphathr{s})}{H} = \norm{H}_\infty$).
Namely, we prove the following concentration result, where, analogously to $t_n$, we let $\ts \in \RR$ denote the quantity such that $2\Phic(\ts) = s/d$.
\begin{proposition}
\label{prop:concentration_gammas}
    There exist universal constants $\cs,\cl,\cd,\cpropf>0$ such that for any $s,d$ with
    $s \geq \cs$ and
    $\cl s \leq d \leq \exp(\cd s^{1/5})$,
\ifarxiv
    \begin{align}
        \abs{
            \frac{\norm{\gamma(\alphathr{s})}_1}{\norm{H}_\infty} - \left(\frac{1}{\ts} - \frac{2}{\ts^3} \right)
        }
        &\leq \frac{\cpropf}{\ts^5}
        & &\text{and} &
        \abs{
            \frac{\norm{\gamma(\alphathr{s})}_2^2}{\norm{H}_\infty^2} - \frac{2}{s \ts^2}
        }
        \leq \frac{\cpropf}{s \ts^4},  \label{eq:concentration}
    \end{align}
\else
    \begin{align}
        \abs{
            \frac{\norm{\gamma(\alphathr{s})}_1}{\norm{H}_\infty} - (\frac{1}{\ts} - \frac{2}{\ts^3})
        }
        &\leq \frac{c}{\ts^5}
        \\
        \abs{
            \frac{\norm{\gamma(\alphathr{s})}_2^2}{\norm{H}_\infty^2} - \frac{2}{s \ts^2}
        }
        &\leq \frac{c}{s t^4}, \label{eq:concentration}
    \end{align}
\fi
    with probability at least $1 - 6\exp \left( -2\frac{s}{\log(d/s)^5} \right)$
    over the draws of $\Hgaussian$.
    %
\end{proposition}
This proposition relies on and extends the literature studying concentration of order statistics \cite{boucheron_2012, li_2020}. An important ingredient for the proof of the proposition is the following lemma, which gives a tight approximation for $\ts$.
\begin{lemma}
\label{lm:tlogds}
    There exist universal constants $\cl, \cone>0$ such that, for all $s \leq d/\cl$, $\ts$ satisfies
    \begin{equation}
        \olt_s^2 - \cone \leq \ts^2 \leq \olt_s^2
    \end{equation}
    where
    \begin{equation}
        \olt_s = \sqrt{2 \log(d/s) - \log \log (d/s) - \log(\frac{\pi}{2})}.
    \end{equation}
    Furthermore, $\cl$ and $\cone$ can be chosen 
    (\eg\ $\cl = 11$ and $\cone = 1$)
    such that
    $\log(d/s) \leq \ts^2 \leq 2 \log(d/s)$.
\end{lemma}
The proofs of Proposition~\ref{prop:concentration_gammas} 
and of Lemma~\ref{lm:tlogds}
are given in Appendix~\ref{apx:subsec:concentration_gamma}.

\subsubsection{Localization: Proof of Proposition~\ref{prop:main_norm} (upper bound for \texorpdfstring{$\phidn$}{phiN})}
\label{subsubsec:proof_bound_phidn}



\label{apx:subsec:bound_phiN}

We now use the concentration bounds of Proposition \ref{prop:concentration_gammas} to obtain a high-probability upper bound for $\phidn$. Recall from Lemma~\ref{lm:parametrization} that it is given by \eqref{eq:phidn_reparam}:
\begin{equation}
    (\phidn)^2 = \min_{\av>0, 1 \leq \alpha \leq \alpham}
    \av^2
    \subjto 
    \av^2 \norm{H}_\infty^2 \geq (1+\cnoise) n (\sigma^2 \norm{\gamma(\alpha)}_1^2 + \av^2 \norm{\gamma(\alpha)}_2^2).
\end{equation}
We may rewrite the constraint as 
\begin{gather}
    \av^2 \norm{H}_\infty^2 \left( 1- (1+\cnoise)n \frac{\norm{\gamma(\alpha)}_2^2)}{\norm{H}_\infty^2} \right)
    \geq (1+\cnoise) n \sigma^2 \norm{\gamma(\alpha)}_1^2 \\
    \iff
    \av^2 \geq 
    \underbrace{
        \frac{\norm{\gamma(\alpha)}_1^2 }{\norm{H}_\infty^2 }
        \frac{
            \sigma^2 n (1+\cnoise)
        }{
            1 - (1+\cnoise) n
            \frac{\norm{\gamma(\alpha)}_2^2}{\norm{H}_\infty^2}
        } 
    }_{=: \ftild(\alpha)^2}
    ~~~\text{and}~~~
    (1+\cnoise) n \frac{\norm{\gamma(\alpha)}_2^2}{\norm{H}_\infty^2} < 1.
\end{gather}
Thus minimizing over $b$ shows that
$(\phidn)^2 = \min_{1 \leq \alpha \leq \alpham} \ftild(\alpha)^2 
\subjto (1+\cnoise) n \frac{\norm{\gamma(\alpha)}_2^2}{\norm{H}_\infty^2} < 1$.
Since we want to upper-bound this minimum, it is sufficient to further restrict the optimization problem by the constraint 
$\alpha \in \left\lbrace \alphathr{s} \middle\vert s \in \{2,...,d\} \right\rbrace$,
yielding
\begin{equation}
\label{eq:proofphilowerbound}
    (\phidn)^2 \leq \min_{2 \leq s \leq d+1}
    \ftild(\alphathr{s})^2
    \subjto (1+\cnoise) n \frac{\norm{\gamma(\alphathr{s})}_2^2}{\norm{H}_\infty^2} < 1.
\end{equation}

We now show that for the choice $s=n$, the constraint is satisfied with high probability, and we give a high-probability estimate for the resulting upper bound $\ftild(\alpha_n)^2$.
See Remark~\ref{rk:choice_s_n} below for a justification of this choice. For the remainder of the proof of Proposition~\ref{prop:main_norm}, we condition on the event where the inequalities in Equation~\eqref{eq:concentration} hold for~$s=n$.
By the concentration bound for $\norm{\gamma(\alphathr{n})}_2^2$, a sufficient condition for the choice $s=n$ to be feasible is
\begin{equation}
    (1+\cnoise)\frac{2}{t_n^2} \left(1+ \frac{c_1}{t_n^2}\right) < 1
\end{equation}
with $c_1>0$ some universal constant.  
Now $t_n^2 \geq \log(d/n)$ by Lemma~\ref{lm:tlogds}, and recall that  $\cnoise = \frac{10}{\log(d/n)^{5/2}}$.
For $\cl$ sufficiently large, the above inequality 
holds for any $n,d$ with $\cl n \leq d$. Moreover, by the concentration bounds for  $\norm{\gamma(\alphathr{n})}_2^2$ and  $\norm{\gamma(\alphathr{n})}_1$,
$f(\alphathr{n})^2$ is upper-bounded by
\begin{align}
  \ftild(\alpha_n)^2
    &\leq
     \frac{
          (1-\frac{4}{t_n^2} + O(\frac{1}{t_n^4}) ) \sigma^2 n (1+\cnoise) 
        }{
            t_n^2 - 2(1+\cnoise) (1+ O(\frac{1}{t_n^2}) )
        }
        \leq
        \frac{\sigma^2 n}{t_n^2} (1+\cnoise) \left(
            1 - \frac{2}{t_n^2} + O(\frac{1}{t_n^4})
        \right).
\end{align}
Furthermore, $\cnoise = O \left( \frac{1}{t_n^5} \right)$ by Lemma~\ref{lm:tlogds}, so 
$f(\alphathr{n})^2 \leq
\frac{\sigma^2 n}{t_n^2} \left(
    1 - \frac{2}{t_n^2} + \frac{\cbrho}{t_n^4}
\right)
=: \Brho(n,d)^2$
for a universal constant $\cbrho>0$.
This concludes the proof of Proposition~\ref{prop:main_norm}.

\begin{remark} \label{rk:choice_s_n}
    Let us informally justify why we can expect the choice $s=n$ to approximately minimize $\ftild(\alpha_s)$. 
    A first justification is that the min-$\ell_1$-norm interpolator $\what$, which is the solution of the optimization problem \eqref{eq:Phi_N}, is well-known to be $n$-sparse. Since the optimization problems \eqref{eq:Phi_N} and \eqref{eq:phidn} are intimately connected via the CGMT (Proposition~\ref{prop:CGMT_application}), we can expect the optimal solution of \eqref{eq:phidn} to have similar properties to $\what$ -- in particular, to have the same sparsity $s=n$. A second, more technical, justification is as follows.
    Note that if we replace 
    $\norm{\gamma(\alphathr{s})}_2^2$
    and
    $\norm{\gamma(\alphathr{s})}_1$
    by their estimates from Proposition~\ref{prop:concentration_gammas}
    and ignore the higher-order terms,
    we have
    \begin{equation}
        \ftild(\alphathr{s})^2
        \approx \frac{1}{t_s^2}
        \frac{
            \sigma^2 n (1+\cnoise)
        }{
            1 - (1+\cnoise) n
            \frac{2}{st_s^2}
        }
        = \frac{
            \sigma^2 (1+\cnoise)
        }{
            \frac{t_s^2}{n} - (1+\cnoise)
            \frac{2}{s}
        }.
    \end{equation}
    Thus a good choice for $s$ is given by maximizing the denominator. By 
    using the estimate $\ts^2 \approx 2 \log(d/s)$ from Lemma~\ref{lm:tlogds},
    we can approximate it by
    \begin{equation}
        \frac{\ts^2}{n} - (1+\cnoise)\frac{2}{s} \approx 
        \frac{2 \log(d/s)}{n} - (1+\cnoise) 
        \frac{2}{s} =: g(s). 
    \end{equation}
    Interpreting $s$ as a continuous variable and setting $\frac{d}{ds}g(s) = 0$ yields the choice $s = (1+\cnoise) n \approx n$.
\end{remark}

\subsubsection{Uniform convergence: Proof of Proposition~\ref{prop:main_upper_lower} (bounds for \texorpdfstring{$\phidp$}{phi+} and \texorpdfstring{$\phidm$}{phi-})}
\label{subsubsec:proof_bound_phidp_phidm}

To obtain an upper bound for the maximization problem defining $\phidp$ (resp. lower bound for the minimization for $\phidm$),
a typical approach would be to find a tractable relaxation of the problem.
However, the more obvious relaxations already explored in the paper \cite{koehler_2021} turn out to be unsatisfactorily loose, as discussed in Section~\ref{sec:proofsketch}.
Here, thanks to the one-dimensional structure of our reformulations \eqref{eq:phidp_reparam} and \eqref{eq:phidm_reparam},
we take a different approach and study the monotonicity of the objectives.

We can decompose our proof in three steps. 
Firstly, we describe our overall monotonicity-based approach.
Secondly, we find values $\ulalpha,\olalpha$ that allow us to unroll our approach.
Finally, we evaluate the bound that the first two steps give us, thus proving the proposition.

\paragraph{Step 1: Studying the feasible set of \eqref{eq:phidp_reparam} and \eqref{eq:phidm_reparam}.}



Recall 
that 
$\phidp$, $\phidm$
are respectively given by 
Equations~\eqref{eq:phidp_reparam},
\eqref{eq:phidm_reparam}.
We can also write them in the following form, using the fact that $\norm{\gamma(\alpha)}_1 = \alpha$:
\begin{align}
    \phidp
    = \max_{1 \leq \alpha \leq \alpham} B^2 \frac{\norm{\gamma(\alpha)}_2^2}{\alpha^2}
    &\subjto 
    B^2 \norm{H}_\infty^2 \geq (1-\cnoise) n (\sigma^2 \alpha^2 + B^2 \norm{\gamma(\alpha)}_2^2) \\
    \phidm
    = \min_{\substack{
        0<\av \leq \gennormboundp \\
        1 \leq \alpha \leq \alpham
    }}
    \av^2 \frac{\norm{\gamma(\alpha)}_2^2}{\alpha^2}
    &\subjto 
    \av^2 \norm{H}_\infty^2 \geq (1-\cnoise) n (\sigma^2 \alpha^2 + \av^2 \norm{\gamma(\alpha)}_2^2).
\end{align}
We first study the sets of feasible solutions of \eqref{eq:phidp_reparam} and \eqref{eq:phidm_reparam}. Denote by $I$ the former set, \ie,
\begin{align}
    I := \left\lbrace
        \alpha \in \left[ 1, \alpham \right]
        ~\middle\vert~
        B^2 \norm{H}_\infty^2 \geq (1-\cnoise) n \left( \sigma^2 \alpha^2 + B^2 \norm{\gamma(\alpha)}_2^2 \right)
    \right\rbrace.
\end{align}
Also let $\alphad = \frac{\norm{H}_1\norm{H}_\infty}{\norm{H}_2^2}$; this choice of notation is purely symbolic, and is justified by the fact that $\alphathr{d} < \alphad < \alphathr{d+1}$.

\begin{lemma} 
\label{lm:interval}
    The following statements hold:
\begin{enumerate}
    \item The mapping $\alpha \mapsto \norm{\gamma(\alpha)}_2^2$ is decreasing over $[1, \alphad]$ and increasing over $[\alphad, \alpham]$.
    \item The mapping $\alpha \mapsto \norm{\gamma(\alpha)}_2^2$ is convex over $\left[1, \alpham \right]$, and $I$ is an interval.
    \item The mapping $\alpha \mapsto \frac{\norm{\gamma(\alpha)}_2^2}{\alpha^2}$ is monotonically decreasing.
\end{enumerate}
\end{lemma}
These monotonicity properties lead us to a proof strategy that can be summarized as follows.

\begin{lemma}
\label{lm:bounds}
    Denote $I = [\ulalpha_I, \olalpha_I]$ the endpoints of $I$.
    For any $\ulalpha \leq \ulalpha_I$,
    \begin{align}
        \phidp &\leq
        \gennormboundp^2 \frac{\norm{\gamma(\ulalpha)}_2^2}{\ulalpha^2}.
    \end{align}
    If $\olalpha_I < \alphad$, then for any $\olalpha_I \leq \olalpha \leq \alphad$,
\ifarxiv
        \begin{align}
            \phidm &\geq \frac{
            \sigma^2 n (1-\cnoise)
        }{
            \norm{H}_\infty^2 - (1-\cnoise) n \norm{\gamma(\oalpha)}_2^2
        } \norm{\gamma(\oalpha)}_2^2. 
    \end{align}
\else
    \begin{align}
        \phidm &\geq \frac{
            \sigma^2 n (1-\cnoise)
        }{
            \norm{H}_\infty^2 - (1-\cnoise) n \norm{\gamma(\oalpha)}_2^2
        } \norm{\gamma(\oalpha)}_2^2. 
    \end{align}
\fi
\end{lemma}
The proofs of Lemma~\ref{lm:interval} and Lemma~\ref{lm:bounds} are given in Appendix~\ref{apx:subsec:feasible_interval}.

\paragraph{Step 2: A tight admissible choice for $\ulalpha$ and $\olalpha$.}

To apply Lemma~\ref{lm:bounds} and obtain bounds on $\phidp$, $\phidm$, we need to find $\ulalpha$ and $\olalpha$ lying on the left, respectively on the right of the interval $I$, and such that $\olalpha \leq \alphad$.
By having a closer look at the way we derived the expression of $\Brho(n,d)$, we have by construction that with high probability, $\alphathr{n} \in I$.
In fact, we show that
there exist integers $\uls$ and $\ols$ very close to $n$
such that
$\alphathr{\uls}$ already falls to the left of $I$, and $\alphathr{\ols}$ to the right of $I$, with high probability.

\begin{lemma}
\label{lm:boundsonalpha}
    Suppose
    ${\norm{\wgt}_1 \leq \cwone \sqrt{ \frac{\sigma^2 n}{\log(d/n)^{5}}}}$ for some universal constant $\cwone>0$.
    There exist universal constants $ \cn, \cl,\cd,\lambda>0$ 
    such that, for any $d,n$ with
    $n \geq \cn$ and
    $\cl n \leq d \leq \exp(\cd n^{1/5})$,
    we can find integers $\us,\os \in \NN_+$ satisfying
\ifarxiv
     \begin{align} \label{eq:boundsons}
         \us =n\exp\left(-\frac{\lambda}{2t_n}\right)\left(1+ O\left(\frac{1}{t_n^2}\right)\right)
        ~~~\text{and}~~~ 
        \os = n \exp\left(\frac{\lambda}{2t_n}\right)\left(1+ O\left(\frac{1}{t_n^2}\right)\right)
    \end{align}
\else
     \begin{align} \label{eq:boundsons}
        \os &= n \exp\left(\frac{\lambda}{2t_n}\right)\left(1+ O\left(\frac{1}{t_n^2}\right)\right) \\
        \text{and}~ 
        \us &=n\exp\left(-\frac{\lambda}{2t_n}\right)\left(1+ O\left(\frac{1}{t_n^2}\right)\right)
    \end{align}
\fi
    and
    \begin{equation} \label{eq:boundsonalpha}
        \alphathr{\us} < \ulalpha_I \leq \alphathr{n} \leq \olalpha_I \leq \alphathr{\os} \leq \alphad,
    \end{equation}
     with probability at least ${1 - 18 \exp\left(-\frac{n}{\log(d/n)^5}\right)}$ over the draws of $\Hgaussian$.
     Moreover, $t_{\us}^2 = t_n^2 + O(1)$ and $t_{\os}^2 = t_n^2 + O(1)$.
\end{lemma}

The proof of the lemma is given in Appendix~\ref{apx:subsec:bound_phip_phim}.
%
It relies in particular on the assumption that 
$\norm{\wgt}_0 \leq \cwone \frac{n}{\log(d/n)^5}$
for some universal constant $\cwone$,
which implies that $\Brho(n,d)$ from Proposition~\ref{prop:main_norm} is the dominating term in $B$, and hence
${
    B^2 = (\Brho(n,d) + c \sigma \sqrt{\sstar})^2
    = \frac{\sigma^2 n}{t_n^2} \left(
        1 - \frac{2}{t_n^2} + O\left( \frac{1}{t_n^4}\right)
    \right)
}$.
Furthermore, the equations in the lemma hold true conditionally on the event where the inequalities in Equation~\eqref{eq:concentration} hold simultaneously for $s=n$, $s=\us$, and $s = \os$ -- which indeed occurs with 
the announced probability.
These two elements of the proof will be reused in the following step.

\paragraph{Step 3: Applying Lemma~\ref{lm:bounds}.}


Lemma~\ref{lm:boundsonalpha} provides us with a choice of $\ulalpha = \alphathr{\uls}$ and $\olalpha = \alphathr{\ols}$ that satisfy the conditions of Lemma~\ref{lm:bounds} with high probability.
To conclude the proof of Proposition~\ref{prop:main_upper_lower}, all that remains to be done is to
compute the bounds given by Lemma~\ref{lm:bounds}, i.e.,
\begin{align}
    \phidp &\leq
    B^2 \frac{\norm{\gamma(\alphathr{\uls})}_2^2}{\norm{\gamma(\alphathr{\uls})}_1^2}
    & &\text{and} &
    \phidm &\geq \frac{
        \sigma^2 n (1-\cnoise)
    }{
        \norm{H}_\infty^2 - (1-\cnoise) n \norm{\gamma(\alphathr{\ols})}_2^2
    } \norm{\gamma(\alphathr{\ols})^2}_2^2. 
\label{eq:upperlowerboundtemp}
\end{align}

For the remainder of the proof of Proposition~\ref{prop:main_upper_lower}, we condition on the event where the inequalities in Equation~\eqref{eq:concentration}
hold simultaneously for $s=n$, $s=\us$, and $s = \os$. 
In particular, the conclusions of Lemma~\ref{lm:boundsonalpha} hold, as discussed just above.
We also recall that, because of the assumption on the growth of $\sstar$, we have
$B^2 = \frac{\sigma^2 n}{t_n^2} \left(
    1 - \frac{2}{t_n^2} + O\left( \frac{1}{t_n^4}\right)
\right)$.
By applying the concentration inequalities from Equation~\eqref{eq:concentration}, and using the above estimate for $B^2$, we obtain
\begin{align}
    \phidp 
    &\leq \frac{\sigma^2 n}{t_n^2} 
    \left(1 + O\left(\frac{1}{t_n^2}\right) \right)
    \frac{2}{\uls t_{\uls}^2 } t_{\uls}^2
    \left(1+O\left(\frac{1}{t_{\uls}^2}\right)\right)
    & &\text{and} &
    \phidm &\geq \frac{2\sigma^2 n}{\ols t_{\ols}^2}
    \left(1 + O\left(\frac{1}{t_n^2}\right)\right). 
\end{align}
By plugging in the approximate expressions of $\uls$ and $\ols$
from Equation~\eqref{eq:boundsons}, 
as well as the estimates $t_{\us}^2 = t_n^2 + O(1)$ and $t_{\os}^2 = t_n^2 + O(1)$ from Lemma~\ref{lm:boundsonalpha},
we further obtain
\begin{align}
    \phidp 
    &\leq \frac{2\sigma^2 }{t_n^2} 
    \exp\left(\frac{\lambda}{2t_n}\right)
    \left(1+ O\left(\frac{1}{t_n^2}\right)\right)  
    & &\text{and}~ &
    \phidm &\geq \frac{2\sigma^2 }{ t_n^2} \exp\left(-\frac{\lambda}{2t_n}\right)
    \left(1 + O\left(\frac{1}{t_n^2}\right)\right). 
\end{align}
Finally, by the expansion $t_n^2 = 2\log(d/n) + O\left(\log\log(d/n)\right)$ from Lemma~\ref{lm:tlogds} and by the Taylor series approximation $\exp(x) = 1+ x + O(x^2)$ (for bounded $x$), we obtain the desired bounds
\begin{align}
    \phidp &\leq \frac{\sigma^2}{\log(d/n)} \left(1+ O\left(\frac{1}{\sqrt{\log(d/n)}}\right)\right)
    & &\text{and} &
    \phidm &\geq \frac{\sigma^2}{\log(d/n)} \left(1+ O\left(\frac{1}{\sqrt{\log(d/n)}}\right)\right).
\end{align}
This concludes the proof of Proposition~\ref{prop:main_upper_lower}.

%% file: sections/conclusion.tex

\ifarx
\section{Future work} 
\else
\section{FUTURE WORK} 
\fi

\label{subsec:future_work}

Our main result gives tight bounds for BP on isotropic Gaussian features.
It would be interesting to extend the study to other connected settings, which we now motivate and for which we summarize key challenges. 
Furthermore, we pose a research question which aims to give a better intuition for the proof. 



\paragraph{Necessity of tightness at the localization step.}
As discussed in Section~\ref{subsec:l1proofsketch}, in order to obtain the right rate in Theorem~\ref{thm:main}, the localization step of our analysis needed to be very tight.
The expression we derive for a high-probability upper bound $M(n,d)$ (from Proposition~\ref{prop:main_norm}) on $\min_{Xw=\xi} \norm{w}_1$ needs to be precise up to relative error of no more than $\Theta \left( \frac{1}{\log(d/n)^2} \right)$.
This strikes us as an unusual feature of our derivation.
Yet, it is unclear whether this is an artifact of our analysis via the application of the GMT, or whether this is due to the nature of the statistical problem itself.
More specifically, we motivate future research to answer the question whether it is true that
\begin{enumerate}
    \item for any $\tilde{M}(n,d) = c M(n,d)$  with $c>1$, we have that with high probability,
    \begin{equation}
        \max_{\substack{
            \norm{w}_1 \leq \tilde{M}(n,d) \\
             Xw=\xi
        }} \norm{w}_2^2 
        \asymp \sigma^2.
    \end{equation}
    \item for any $\tilde{M}(n,d) = M(n,d)(1+\omega(\frac{1}{\log(d/n)}))$, we have that with high probability,
    \begin{equation}
    \max_{\substack{
        \norm{w}_1 \leq \tilde{M}(n,d) \\
         Xw=\xi
    }} \norm{w}_2^2 
    = \omega\left(\frac{\sigma^2}{\log(d/n)}\right).
    \end{equation}
\end{enumerate}


Resolving this question is challenging due to the non-concavity of the maximization objective. While we can still use the GMT to upper-bound this quantity (see Proposition~\ref{prop:CGMT_application}), we cannot use the CGMT to lower-bound it, and thus the methodologies
used in this paper fall short.
As a possible direction, we note that
this hypothesis is related to the question of finding tight lower bounds for the diameter of the intersection of the kernel of $X$ and the $\ell_1$-ball (see Theorem~3.5 in \cite{vershynin2011lectures}).

\paragraph{Non-isotropic features.}
Theorem~\ref{thm:main} 
assumes isotropic features as we are interested in
showing consistency of BP for inherently high-dimensional input data.
By contrast, recently there has been an increased interest in studying spiked covariance data models (see e.g.\ \cite{bartlett_2020,muthukumar2021classification,chatterji_2021}). In such settings even  
min-$\ell_2$-norm interpolators can achieve consistency. 
The main obstacle 
to extending our methodology to 
non \iid\ features lies in 
adapting the definition 
of the path $\{\gamma(\alpha)\}_\alpha$.
Assuming a diagonal covariance matrix,
such an extension should be relatively straightforward.
We leave this task and the challenging non-diagonal case for future work. 


\paragraph{Non-Gaussian features.}

The proof of Theorem~\ref{thm:main} crucially relies on the (Convex)
Gaussian Minimax Theorem \cite{thrampoulidis_2015,gordon_1988}, and
hence on the assumption that the input features are drawn from a
Gaussian distribution. In Figure~\ref{fig:distributions}, we include
plots of the prediction error $\norm{\what - \wgt}_2^2$ not only for Gaussian but also for Log Normal and Rademacher distributed features.  We
observe that in all three cases, the prediction error 
closely follows the trend line $\frac{\sigma^2}{\log(d/n)}$ (dashed curve). This leads 
us to
conjecture that Theorem~\ref{thm:main} can be extended to a more general class of
distributions.

Generalizing our results in this direction appears to be a
challenging task since
the tools 
used
in this paper are not directly applicable anymore.
Instead, for heavy-tailed distributions, a popular theoretical framework is the small-ball method
\cite{mendelson_2014,koltchinskii_2015}, which covers the Log Normal and Rademacher distributions. 
\ifarx The authors of the paper
\cite{chinot_2021} 
\else
\citet{chinot_2021}
\fi
apply this approach to
min-$\ell_1$-norm interpolation, and obtain the constant upper bound $O(\sigma^2)$, under more general assumptions than our setting 
(in particular their analysis handles adversarial noise with magnitude controlled by $\sigma^2$).
Yet, it is unclear whether the looseness of their upper bound is an artifact of their proof, or whether the small-ball method itself is too general to capture the rates observed in Figure~\ref{fig:distributions}.

Finally,
we also leave
it as future work to adapt our proof technique for minimum-norm interpolators with general norms,
and for classification tasks.



\ifarx
\section{Conclusion}
\else
\section{CONCLUSION}
\fi
\label{sec:conclusion}
By introducing a novel proof technique, we derive matching upper and lower bounds of order $\frac{\sigma^2}{\log(d/n)}$ on the prediction error of 
basis pursuit (BP, or min-$\ell_1$-norm interpolation) in noisy sparse linear regression.
Our result closes a gap in the minimum-norm interpolation literature,
disproves a conjecture from \cite{chinot_2021},
and is the first to imply asymptotic consistency of a minimum-norm interpolator
for isotropic features.
Furthermore, the prediction error decays with the amount of overparametrization $d/n$, confirming that BP also benefits from the regularization effect of high dimensionality, as suggested by the modern storyline on interpolating models.

%% file: sections_proof/l1_proof.tex

\ifarxiv
\else
\onecolumn
\fi

\ifarxiv
\section{Proof details}
\else
\section{PROOF DETAILS}
\fi
\label{apx:sec:appendix_proofs}



In this appendix, we provide details of the proof of our main result, Theorem~\ref{thm:main}, omitted in Section~\ref{sec:proof_maintext}.
We refer to that section for notation.



\subfile{0_prooftriangle}

\subfile{1_CGMT}

\subfile{2_param_argmaxmin_v2}

\subfile{3_geometric_lemma}

\subfile{4_concentration_gamma}


\subfile{6_feasible_interval}

\subfile{7_bound_phip_phim}

%% file: sections_proof/0_prooftriangle.tex

\subsection{Proof of Lemma~\ref{lm:triangle_ineq}: Preliminary}
\label{apx:subsec:triangle}

Let $\Sparseset := \{j:~ \wgt_j \neq 0\}$
(in particular, $\sstar = \vert \Sparseset \vert$).
Denote by $\whatS \in \RR^d$ the vector with entries 
$(\whatS)_j = \what_j$ if $j \in \Sparseset$ and $0$ otherwise,
and let $\whatSm = \what - \whatS$. 
Now recall the definition of $\what$,
and define $\nuhat$ by
\begin{align}
    \what &= \arg\min_w \norm{w}_1 \subjto X(w-\wgt) = \xi \\
    \nuhat &= \arg\min_v \norm{v}_1 \subjto Xv = \xi.
\end{align}

Clearly by definition of $\nuhat$ and $\what$, we have that ${\norm{\what}_1 \leq \norm{\nuhat}_1 + \norm{\wgt}_1}$.
Therefore
\begin{align}
    0 &\leq \norm{\wgt}_1 - \norm{\what}_1 + \norm{\nuhat}_1\\
    &= \norm{\wgt}_1 - \norm{\whatS}_1 - \norm{\whatSm} + \norm{\nuhat}_1\\
    &\leq \norm{\wgt - \whatS}_1 - \norm{\whatSm}_1 + \norm{\nuhat}_1\\
    &= 2\norm{\wgt - \whatS}_1 - \norm{\what - \wgt }_1  + \norm{\nuhat}_1\\
    &\leq 2\sqrt{\sstar} \norm{\wgt - \whatS}_2  - \norm{\what - \wgt }_1  + \norm{\nuhat}_1.
\end{align}
Hence, 
\begin{align}
    \norm{\what - \wgt }_1
    &\leq 2\sqrt{\sstar} \norm{\wgt - \whatS}_2  + \norm{\nuhat}_1 \\
    &\leq 2\sqrt{\sstar} \norm{\wgt - \what}_2  + \norm{\nuhat}_1.
\end{align}
Finally, we bound $\norm{\wgt - \what}_2$
by applying Theorem~3.1 of \cite{chinot_2021}, noting that its assumptions are subsumed by the assumptions of our Theorem~\ref{thm:main}.

%% file: sections_proof/1_CGMT.tex

\subsection{Proof of Proposition~\ref{prop:CGMT_application}: Application of the (C)GMT}
\label{apx:subsec:CGMT}

Proposition~\ref{prop:CGMT_application} reduces the estimation of the quantities $\Phi_N, \Phi_+, \Phi_-$ in Equations~\eqref{eq:Phi_N}, \eqref{eq:Phi_+}, \eqref{eq:Phi_-} to the estimation of auxiliary quantities $\phidn, \phidp, \phidm$, using the (C)GMT.

As a first step, we apply the CGMT to $\Phi_N$ and the GMT to $\Phi_+$ analogously to \cite[Lemmas~4\&7]{koehler_2021}. We only restate the results here and refer the reader to that paper for details and proofs.
Note that the (C)GMT is applied on $X$ conditionally on $\xi$, so that the Gaussianity of the noise is not crucial.

\begin{lemma}[{\cite[Lemma~7]{koehler_2021}, Application of CGMT}]
    Define
    \begin{equation}
        \tphi_N = \min_w \norm{w}_1 \subjto 
        \norm{\xi - \Ggaussian \norm{w}_2}_2 \leq \innerprod{w}{\Hgaussian},
    \end{equation}
    where $\Ggaussian \sim \NNN(0,I_n)$ and $\Hgaussian \sim \NNN(0,I_d)$ are independent random variables.
    Then, for all $t \in \RR$,
    \begin{equation}
        \PP(\Phi_N>t) \leq 2 \PP(\tphi_N >t),
    \end{equation}
    where the probabilities
    on the left and on the right
    are over the draws of $X,\xi$ and of $\Ggaussian,\Hgaussian,\xi$, respectively. 
\end{lemma}

\begin{lemma}[{\cite[Lemma~4]{koehler_2021}, Application of GMT}]
    Define
    \begin{equation}
        \tphi_+ = \max_w \norm{w}_2^2
        \subjto \begin{cases}
            \norm{w}_1  \leq B(n,d) \\
            \norm{\xi - \Ggaussian \norm{w}_2}_2 \leq \innerprod{w}{\Hgaussian}
        \end{cases},
    \end{equation}
    where $\Ggaussian \sim \NNN(0,I_n)$ and $\Hgaussian \sim \NNN(0,I_d)$ are independent random variables.
    Then, for all $t \in \RR$,
    \begin{equation}
        \PP(\Phi_+>t) \leq 2 \PP(\tphi_+ >t),
    \end{equation}
    where the probabilities
    on the left and on the right
    are over the draws of $X,\xi$ and of $\Ggaussian,\Hgaussian,\xi$, respectively. 
\end{lemma}

Following the same argument as in \cite[Lemma~4]{koehler_2021}, we can show a corresponding lemma for $\Phi_-$ which we state without proof:
\begin{lemma}[{Application of GMT}]
    Define
    \begin{equation}
        \tphi_- = \min_w \norm{w}_2^2
        \subjto \begin{cases}
            \norm{w}_1  \leq B(n,d) \\
            \norm{\xi - \Ggaussian \norm{w}_2}_2 \leq \innerprod{w}{\Hgaussian}
        \end{cases},
    \end{equation}
    where $\Ggaussian \sim \NNN(0,I_n)$ and $\Hgaussian \sim \NNN(0,I_d)$ are independent random variables.
    Then, for all $t \in \RR$,
    \begin{equation}
        \PP(\Phi_-<t) \leq 2 \PP(\tphi_- <t),
    \end{equation}
    where the probabilities
    on the left and on the right
    are over the draws of $X,\xi$ and of $\Ggaussian,\Hgaussian,\xi$, respectively. 
\end{lemma}

Next, by using Gaussian concentration results, we can formulate simpler versions of the above optimization problems defining the quantities $\tphi$'s.

\paragraph{Simplify $\tphi_N$.}

Following the same argument as in the first part of the proof of \cite[Lemma~8]{koehler_2021} (Equations~(68)-(70)),
we can show that for any $0<\rho<1/2$, with probability at least $1-6\exp(-n \cnoise^2/100)$,
uniformly over $w$,
\begin{equation}
    \norm{\xi - \Ggaussian \norm{w}_2}_2^2
    \leq 
    (1+\cnoise) n (\sigma^2 + \norm{w}_2^2).
\label{eq:apx1eventUpper}
\end{equation}
So on the event where Equation~\eqref{eq:apx1eventUpper} holds, we have that
\begin{equation}
    \tphi_N \leq
        \min_w \norm{w}_1 \subjto 
        \innerprod{w}{\Hgaussian}^2 \geq (1+\cnoise) n (\sigma^2 + \norm{w}_2^2)
    = \phidn
\end{equation}
%
which proves the
first inequality
in  Proposition~\ref{prop:CGMT_application}.

\paragraph{Simplify $\tphi_+,\tphi_-$.}

By the same argument as for $\tphi_N$,
we can show that for any $0<\rho<1/2$, with probability at least $1-6\exp(-n \cnoise^2/100)$,
uniformly over $w$,
\begin{equation}
    \norm{\xi - \Ggaussian \norm{w}_2}_2^2
    \geq 
    (1-\cnoise) n (\sigma^2 + \norm{w}_2^2).
\label{eq:apx1eventLower}
\end{equation}
So on the event where Equation~\eqref{eq:apx1eventLower} holds, we have that
\begin{equation}
    \tphi_+ \leq
    \max_w \norm{w}_2^2
    \subjto \begin{cases}
        \norm{w}_1  \leq M+2\|\wgt\|_1 \\
        \innerprod{w}{\Hgaussian}^2 \geq (1-\cnoise) n (\sigma^2 + \norm{w}_2^2)
    \end{cases} 
    = \phidp,
\end{equation}
and similarly $\tphi_- \geq \phidm$.
This proves the
second and third inequalities
in Proposition~\ref{prop:CGMT_application} and thus completes the proof.

%% file: sections_proof/2_param_argmaxmin_v2.tex
\subsection{Proof of Lemma~\ref{lm:parametrization}: Parametrizing the argmax/argmin}
\label{apx:subsec:param_argmaxmin}

We now prove our first key lemma: we show that, up to scaling, the argmax/argmin in 
\eqref{eq:phidn},
\eqref{eq:phidp}, and
\eqref{eq:phidm}
belong to a certain parametric path $\Gamma = \{ \frac{\gamma(\alpha)}{\alpha} \}_\alpha$.
Throughout this section and the next, we consider $H$ as a fixed vector such that
\mbox{$H_1 > ... > H_d > 0$}.
In other words, all of our statements should be understood as holding conditionally on $\Hgaussian$, and with $\Hgaussian$ in general position.

For all $\beta \in [\frac{1}{d}, 1]$, define
\begin{equation} \label{eq:def_olw_beta}
    \olw(\beta) = \argmax_w \innerprod{w}{H} \subjto \begin{cases}
        \norm{w}_2^2 \leq \beta \\
        \forall i, w_i \geq 0 \\
        \onevec^\top w = \norm{w}_1 = 1
    \end{cases}.
\end{equation}
Importantly, note that the constraint 
$\norm{w}_2^2 \leq \beta$
in the definition of $\olw(\beta)$ necessarily holds with equality at optimality.
Indeed, suppose by contradiction 
$\norm{\olw(\beta)}_2^2 < \beta \leq 1 = \norm{\olw(\beta)}_1^2$.
This implies that $\olw(\beta)$ has at least two nonzero components; denote $i \neq 1$ such that $\olw(\beta)_i > 0$.
Then there exists some $\eps>0$ such that $\olw(\beta) + \eps e_1 - \eps e_i$ satisfies the constraints and achieves a higher objective value than $\olw(\beta)$, contradicting its optimality.

The first step of the proof is to show that (at least one element) of the argmax/argmin belongs to the set $\RR_+ \olGamma$,
where $\olGamma = \left\lbrace
    \olw(\beta) ; \frac{1}{d} \leq \beta \leq 1
\right\rbrace$.

\begin{claim}
    For each of the optimization problems 
    \eqref{eq:phidn},
    \eqref{eq:phidp}, and
    \eqref{eq:phidm},
    there exist $b>0$ and
    $\beta \in [\frac{1}{d}, 1]$
    such that $b \olw(\beta)$ is an optimal solution.

    
    
\end{claim}

\begin{proof}
    Let $v$ be an optimal solution of \eqref{eq:phidn}.
    It is straightforward to check that we may assume \WLOG\ that $\forall i, v_i \geq 0$.
    Choose $b = \norm{v}_1$
    and $\beta = \frac{\norm{v}_2^2}{\norm{v}_1^2}$;
    note that $\beta \in [\frac{1}{d}, 1]$.
    By definition,
    $\norm{b \olw(\beta)}_2 = \norm{v}_2$ and
    $\norm{b \olw(\beta)}_1 = \norm{v}_1$,
    and $v/b$ is feasible for \eqref{eq:def_olw_beta} so
    $\innerprod{b \olw(\beta)}{H} \geq \innerprod{v}{H}$.
    Therefore, $b \olw(\beta)$ satisfies the constraint of \eqref{eq:phidn} and achieves the optimal objective value, so is also an optimal solution of \eqref{eq:phidn}. 
    
    The statements for \eqref{eq:phidp} and \eqref{eq:phidm} follow by the exact same argument.
\end{proof}

Next, we show that $\{ \olw(\beta) \}_\beta$ and $\{ \frac{\gamma(\alpha)}{\alpha} \}_\alpha$ are two parametrizations of the same path. 

\begin{claim}
    We have the equality
    \begin{equation}
        \olGamma :=
        \left\lbrace \olw(\beta) ; \beta \in \left[\frac{1}{d}, 1\right] \right\rbrace
        ~~=~~
        \left\lbrace \frac{\gamma(\alpha)}{\alpha} ; \alpha \in [1, 
        \alpham] \right\rbrace
        =: \Gamma
    \end{equation}
    where $\alpham = d \frac{\norm{H}_\infty}{\norm{H}_1}$.
\end{claim}

\begin{proof}
    First note that we can characterize $\frac{\gamma(\alpha)}{\alpha}$ as the optimal solution of
    \begin{equation} \label{eq:def_gammap_alphap}
        \frac{\gamma(\alpha)}{\alpha} = \argmin_w \norm{w}_2^2
        \subjto \begin{cases}
            \innerprod{w}{H} \geq \frac{\norm{H}_\infty}{\alpha} \\
            \forall i, w_i \geq 0 \\
            \onevec^\top w = \norm{w}_1 = 1
        \end{cases}.
    \end{equation}
    
    The optimization problems \eqref{eq:def_olw_beta} and \eqref{eq:def_gammap_alphap} are both convex and both satisfy the Linear Independence Constraint Qualification conditions.
    So, denoting $\Delta_d = \left\lbrace
        w \in [0,1]^d ; \onevec^\top w =1
    \right\rbrace$
    the standard simplex,
    by the Lagrangian duality theorem (a.k.a.\ Karush-Kuhn-Tucker theorem) we have that for all $w \in \RR^d$,
    \begin{align}
        &\exists \beta>0; w = \olw(\beta)
        &
        &\iff 
        &
        &\exists \lambda>0; w = \argmax_{w \in \Delta_d} \innerprod{w}{H} - \lambda \norm{w}_2^2 &&&& \\
        && 
        &\iff 
        &
        &\exists \mu>0; w = \argmin_{w \in \Delta_d} \norm{w}_2^2 - \mu \innerprod{w}{H}
        &
        &\iff
        &
        &\exists \alpha>0; w = \frac{\gamma(\alpha)}{\alpha}.
    \end{align}
    
    Thus, 
    ${
        \left\lbrace \olw(\beta) ; \beta>0 \right\rbrace
        =
        \left\lbrace \frac{\gamma(\alpha)}{\alpha} ; \alpha>0 \right\rbrace
    }$.
    However it is straightforward to check that
    ${
        \left\lbrace \olw(\beta) ; \beta>0 \right\rbrace
        = 
        \olGamma
    }$
    and that
    ${
        \left\lbrace \frac{\gamma(\alpha)}{\alpha} ; \alpha>0 \right\rbrace
        = 
        \Gamma
    }$,
    which concludes the proof.
\end{proof}

Just as the first constraint in \eqref{eq:def_olw_beta} holds with equality at optimality, so does the first constraint in \eqref{eq:def_gammap_alphap}; that is,
$\innerprod{\gamma(\alpha)}{H} = \norm{H}_\infty$ for all $\alpha \in [1, \alpham]$.
This would follow from a careful study of the equivalence between the two problems, but here we give a more direct proof.

\begin{claim} \label{claim:saturate_in_gamma}
    The inequality constraint
    $\innerprod{w}{H} \geq \norm{H}_\infty$
    in the problem defining $\gamma(\alpha)$ holds with equality, at optimality.
\end{claim}

\begin{proof}
    Denote $w = \gamma(\alpha)$.
    Suppose by contradiction 
    $\innerprod{w}{H} > \norm{H}_\infty$.
    Let $i$ resp.\ $j$ the index of the largest resp.\ smallest component of $w$.
    First note that if $w_i = w_j$, then $w \propto \bmone$, \ie, $w = \frac{\alpha}{d} \bmone$ and so $\innerprod{w}{H} = \frac{\alpha}{d} \norm{H}_1 > \norm{H}_\infty$, which would contradict $\alpha \leq \alpham$;
    so we have the strict inequality $w_i > w_j$.
    Now for some $\eps>0$ to be chosen, let $w' = w - \eps e_i + \eps e_j$.
    Clearly $\eps>0$ can be chosen small enough so that $w'$ satisfies all three constraints in the optimization problem defining $\gamma(\alpha)$.
    Furthermore, for small enough $\eps$,
    $\norm{w}_2^2 - \norm{w'}_2^2 = w_i^2 - (w_i-\eps)^2 + w_j^2 - (w_j+\eps)^2 = 2\eps \left( w_i-w_j - \eps \right)$ is positive, i.e $\norm{w'}_2^2 < \norm{w}_2^2$, which contradicts optimality of $w = \gamma(\alpha)$.
\end{proof}

We now have all the necessary ingredients to prove Lemma~\ref{lm:parametrization}.
The equivalence between \eqref{eq:phidn} and \eqref{eq:phidn_reparam} follows immediately from constraining the variable $w$ (in the former) to belong to the set $\RR_+ \Gamma$.
The equivalence between \eqref{eq:phidm} and \eqref{eq:phidm_reparam} also follows immediately, noting that $\norm{b \frac{\gamma(\alpha)}{\alpha}}_1 = b$ by definition.
Finally, the equivalence between \eqref{eq:phidp} and \eqref{eq:phidp_reparam} follows by noticing that the inequality constraint $\norm{w}_1 \leq B$ (in the former) is necessarily saturated at optimality.

%% file: sections_proof/3_geometric_lemma.tex
\subsection{Proof of Lemma~\ref{lm:charact_gamma}: Characterizing the parametric path}
\label{apx:subsec:geometric_lemma}

We now
give a precise characterization of the parametric path $\Gamma$, by studying the optimization problem defining $\gamma(\alpha)$.
%
%
%
Throughout this section (just as in the previous one), we consider $H$ as a fixed vector such that
\mbox{$H_1 > ... > H_d > 0$}.
In other words, all of our statements should be understood as holding conditionally on $\Hgaussian$, and with $\Hgaussian$ in general position.





Throughout the proof, consider a fixed $1<\alpha \leq \alpham$. 
The goal is to derive a closed-form expression of $\gamma(\alpha)$.
We proceed by a Lagrangian duality approach, and first identify the dual variables $\lambda, \nu, \mu$ (a.k.a.\ Lagrangian multipliers, a.k.a.\ KKT vectors) of the optimization problem defining $\gamma(\alpha)$.
This first analysis yields an expression for $\gamma(\alpha)$ involving an unknown ``sparsity'' integer $s$, which depends on $\lambda, \nu, \mu$ and hence indirectly on $\alpha$.
We finish by showing how to determine $s$ explicitly from $\alpha$.

\paragraph{Karush-Kuhn-Tucker (KKT) conditions.}

Recall that $\gamma(\alpha)$ is defined by the following optimization problem (note that the additional factor $\frac{1}{2}$ in the objective does not change the $\argmin$):
\begin{equation}
    \gamma(\alpha) = \argmin_w \frac{1}{2} \norm{w}_2^2
    \subjto \begin{cases}
        \innerprod{w}{H} \geq \norm{H}_\infty \\
        \forall i, w_i \geq 0 \\
        \onevec^\top w = \norm{w}_1 = \alpha
    \end{cases}.
\end{equation}
This is a convex optimization problem with Lagrangian
\begin{equation}
    L(w; \lambda, \mu, \nu) = \frac{1}{2} \norm{w}_2^2 
    - \lambda (\innerprod{w}{H} - \norm{H}_\infty)
    + \mu (\onevec^\top w - \alpha)
    - \nu^\top w.
\end{equation}
The objective is convex and all the constraints are affine.
So by Lagrangian duality, $w = \gamma(\alpha)$ if and only if there exist
$\lambda, \mu \in \RR$ and $\nu \in \RR^d$
satisfying the KKT conditions:
\begin{itemize}
    \item (Stationarity)
    \mbox{$w - \lambda H + \mu \onevec - \nu = 0$}
    ~\ie~
    \mbox{$w = \lambda H - \mu \onevec + \nu$}
    \item (Primal feasibility)
    \mbox{$\begin{cases}
        \innerprod{w}{H} \geq \norm{H}_\infty \\
        \forall i,~ w_i \geq 0 \\
        \onevec^\top w = \alpha
    \end{cases}$}
    \item (Dual feasibility) 
    \mbox{$\lambda \geq 0$}
    and
    \mbox{$\forall i,~ \nu_i \geq 0$}
    \item (Complementary slackness)
    \mbox{$\lambda (\innerprod{w}{H} - \norm{H}_\infty) = 0$},
    and
    \mbox{$\forall i,~ \nu_i w_i = 0$}.
\end{itemize}
In the rest of this proof, denote $w = \gamma(\alpha)$, and let $\lambda, \mu, \nu$ as above.

\paragraph{Sparsity structure of $w$.}

Let $s$ denote the largest $s' \in \{1, ..., d \}$ such that $\lambda H_{s'} > \mu$.
Since $\lambda \geq 0$ and $H$ is ordered, we have 
\begin{equation}
    \lambda H_1 \geq ... \geq \lambda H_s > \mu \geq \lambda H_{s+1} \geq ...
\end{equation}
Consider the complementary slackness condition $\forall i,~ \nu_i w_i = 0$.
\begin{itemize}
    \item If $w_i > 0$, then $\nu_i=0$ so $w_i = \lambda H_i - \mu > 0$, and so $i \leq s$. 
    \item If $\nu_i > 0$, then $w_i = \lambda H_i - \mu + \nu_i = 0$ so $\lambda H_i - \mu < 0$, and so $i > s$.
    \\
    So by contraposition, for all $i \leq s$, $\nu_i=0$ and $w_i = \lambda H_i - \mu$.
\end{itemize}
Thus, $\support(w) \subset \{1, ..., s\}$ and $w = \lambda \Hvecs - \mu \onevecs$,
where $\Hvecs$ is the vector equal to $H$ on the first $s$ components and $0$ on the last $(d-s)$, and similarly for $\onevecs$.

Furthermore, note that the case $s=1$ occurs only if $w \propto e_1$, and one can check that it implies $\alpha = 1$, which we excluded.

\paragraph{Closed-form expression of the dual variables $\lambda,\mu$.}

We can compute $\lambda$ and $\mu$ by substituting $w = \lambda \Hvecs - \mu \onevecs$ into the primal feasibility conditions.
\begin{itemize}
    \item Since we know from Claim~\ref{claim:saturate_in_gamma}
    (in Section~\ref{apx:subsec:param_argmaxmin}) that the first constraint in the problem defining $\gamma(\alpha)$ holds with equality at optimality, this means that the first primal feasibility condition holds with equality, \ie,
    \begin{equation}
        \innerprod{w}{H} = \lambda \norm{\Hvecs}_2^2 - \mu \norm{\Hvecs}
        =
        \norm{H}_\infty.
    \end{equation}
    \item By the last primal feasibility condition,
    \mbox{$\onevec^\top w = \lambda \norm{\Hvecs}_1 - \mu s 
    =
    \alpha$}.
\end{itemize}
So $\lambda$ and $\mu$ are given by
\begin{equation}
    \begin{cases}
        \lambda \norm{\Hvecs}_2^2 - \mu \norm{\Hvecs}
        &=
        \norm{H}_\infty \\
        \lambda \norm{\Hvecs}_1 - \mu s
        &=
        \alpha
    \end{cases}
    \iff
    \begin{cases}
        \lambda &= \frac{1}{s \norm{\Hvecs}_2^2 - \norm{\Hvecs}_1^2}
        \left( s \norm{H}_\infty - \alpha \norm{\Hvecs}_1 \right) \\
        \mu &= \frac{1}{s \norm{\Hvecs}_2^2 - \norm{\Hvecs}_1^2}
        \left( \norm{\Hvecs}_1 \norm{H}_\infty - \alpha \norm{\Hvecs}_2^2 \right)
    \end{cases}.
\end{equation}
Note that the denominator is positive, since $\Hvecs$ has distinct components.

\paragraph{Closed-form characterization of $s$.}
We now show that there exists an increasing sequence $\alpha_2 = 1 < \cdots < \alpha_d <\alpha_{d+1} = \alpham $ such that for all $\alpha$, the sparsity  $s$ of $w$ is exactly the index which satisfies $\alpha \in (\alpha_s,\alpha_{s+1}] $.
By plugging the expressions of $\lambda$ and $\mu$ into the condition defining $s$:
$\lambda H_s > \mu \geq \lambda H_{s+1}$,
we obtain
\begin{equation}
    \left( s \norm{H}_\infty - \alpha \norm{\Hvecs}_1 \right) H_s
    >
    \norm{\Hvecs}_1 \norm{H}_\infty - \alpha \norm{\Hvecs}_2^2
    \geq
    \left( s \norm{H}_\infty - \alpha \norm{\Hvecs}_1 \right) H_{s+1}.
\end{equation}
Rearranging, this is equivalent to
\begin{align}
    \alpha \left( \norm{\Hvecs}_2^2 - \norm{\Hvecs}_1 H_s \right)
    &>
    \left( \norm{\Hvecs}_1 - s H_s \right) \norm{H}_\infty \\
    \text{and}~~
    \alpha \left( \norm{\Hvecs}_2^2 - \norm{\Hvecs}_1 H_{s+1} \right)
    &\leq
    \left( \norm{\Hvecs}_1 - s H_{s+1} \right) \norm{H}_\infty.
\end{align}
One can check that 
\mbox{$
\norm{\Hvecs}_2^2 - \norm{\Hvecs}_1 H_{s+1}
>
\norm{\Hvecs}_2^2 - \norm{\Hvecs}_1 H_s
> 0
$}.
So the above is equivalent to
\begin{equation}
    \alphathr{s} := 
    \frac{
        \left( \norm{\Hvecs}_1 - s H_s \right) \norm{H}_\infty
    }{
        \norm{\Hvecs}_2^2 - \norm{\Hvecs}_1 H_s
    }
    < \alpha \leq
    \frac{
        \left( \norm{\Hvecs}_1 - s H_{s+1} \right) \norm{H}_\infty
    }{
        \norm{\Hvecs}_2^2 - \norm{\Hvecs}_1 H_{s+1}
    }
    =: \olalpha(s).
\end{equation}
A straightforward calculation shows that
$\olalpha(s) = \alphathr{s+1}$.
Thus, using the convention 
$\alphathr{d+1} = \alpham$,
$s$ is uniquely characterized by
$\alphathr{s} < \alpha \leq \alphathr{s+1}$.
This concludes the proof of Lemma~\ref{lm:charact_gamma}.

%% file: sections_proof/4_concentration_gamma.tex
\subsection{Proof of Proposition~\ref{prop:concentration_gammas}: Concentration of norms of  \texorpdfstring{$\gamma(\alphathr{s})$}{gamma(alphas)}}
\label{apx:subsec:concentration_gamma}
In this section we prove Proposition~\ref{prop:concentration_gammas} and Lemma~\ref{lm:tlogds} which establish concentration inequalities for $\gamma(\alpha)$ at the breakpoints
$\alphathr{s}$ for $2 \leq s \leq d$.
More precisely, we give high-probability estimates (with respect to the draws of $\Hgaussian$) of their $\ell_1$ and $\ell_2$ norms, since those are the quantities that appear in
the stochastic optimization problems
\eqref{eq:phidn_reparam},
\eqref{eq:phidp_reparam} and
\eqref{eq:phidm_reparam}.

Plugging in $\alpha = \alphathr{s}$ into the closed-form expressions of $\lambda$ and $\mu$ in Lemma~\ref{lm:charact_gamma}, we obtain
${\gamma(\alphathr{s}) = \frac{\norm{H}_\infty}{\langle v_s, H\rangle} v_s}$
where $v_s := H_{[s]}- H_s \onevecs $.
Thus, to estimate the norms of $\gamma(\alphathr{s})$ it suffices to estimate the quantities
\begin{align}
    \norm{v_s}_2^2 &= \norm{H_{[s]}}_2^2 - 2 \norm{H_{[s]}}_1 H_s + s H_s^2 \\
    \norm{v_s}_1 &= \norm{H_{[s]}}_1 - s H_s \\
    \innerprod{v_s}{H} &= \norm{H_{[s]}}_2^2 -  \norm{H_{[s]}}_1 H_s.
\end{align}

Throughout the proofs in this section, we will use $\gencstpos>0$ to denote a universal constant (in particular, independent of $d$ and $s$) which may change from display to display.
Furthermore, in this section we let $Z$ denote a standard normal distributed random variable,
and recall that $\Phic(x) = \PP(Z>x) = \frac{1}{2} \PP(\abs{Z}>x)$ for $x>0$ denotes its complementary cumulative distribution function.

\subsubsection{Preliminary facts}

We start by stating some auxiliary facts about $\Phic$.
\begin{fact} \label{fact:estim_Phic}
    Denote $h(x)$ the function such that 
    $\forall x>0,~ \Phic(x) = \frac{\exp(-x^2/2)}{x \sqrt{2\pi}} h(x)$.
    We have the first-order and higher-order upper and lower bounds
    \begin{align}
        1-\frac{1}{1+x^2} &\leq h(x) \leq 1
        ~~~\mathrm{and}~~~
        \abs{
            h(x) - \left( 1- \frac{1}{x^2} + \frac{3}{x^4}  - \frac{15}{x^6} \right)
        }
        \leq \frac{c}{x^8}
    \end{align}
    for all $x>0$.
\end{fact}
\begin{proof}
    The first-order estimate follows from straightforward analysis.
    The higher-order estimate follows from the exact asymptotic expansion of the \textit{complementary error function} $\text{erfc}$,%
    since
    $2\Phic(x) = \text{erfc}(x/\sqrt{2})$.
\end{proof}

\begin{fact} \label{fact:cond_expectations}
    By straightforward calculations, we have
    \begin{align}
        \forall x>0,~
        \EE [Z \vert Z \geq x]
        &= \frac{1}{\Phic(x)} \frac{\exp(-x^2/2)}{\sqrt{2\pi}}
        = \frac{x}{h(x)} \\
        \mathrm{and}~~~
        \EE[Z^2 | Z \geq x] 
        &= \frac{1}{\Phic(x)} \left(
            \frac{x}{\sqrt{2\pi}} \exp(-x^2/2) + \Phic(x)
        \right)
        = 1+\frac{x^2}{h(x)}.
    \end{align}
\end{fact}

We will also make repeated use of Lemma~\ref{lm:tlogds} (Section~\ref{subsubsec:concentration_gammas}), whose proof is deferred to Section~\ref{apx:subsubsec:concentration_gamma:tlogs}.

\subsubsection{Proof of Proposition~\ref{prop:concentration_gammas}}

We will show the following lemmas successively, in which $t \in \RR$ denotes the quantity such that ${2\Phic(t) = s/d}$
(we drop the explicit dependency on $s$ for concision in this section).

\begin{lemma}[Concentration of $H_s$] \label{lm:concentration_Hs}
    Assume that $s < d/2$.
    With probability at least $1-2\delta$, we have
    \begin{equation}
        \abs{H_s - t} \leq \gencstpos \left( 
            \frac{1}{\sqrt{s}} 
            + \sqrt{\frac{\log(1/\delta)}{s}}
            + \frac{\log(1/\delta)}{s}
        \right).
    \end{equation}
\end{lemma}

\begin{lemma}[Concentration of $\norm{H_{[s]}}_2^2$] \label{lm:concentration_l2_Hs}
    Assume $s < d/5$. 
    With probability at least $1-2\delta$, we have
    \begin{equation}
        \abs{
            \norm{\Hvecs}_2^2 - s \EE[Z^2|Z \geq t]
        } \leq 
        \gencstpos \sqrt{s} (1+\sqrt{\log(1/\delta)})
        \left(
            \frac{1}{\sqrt{s}} (1+\sqrt{\log(1/\delta)})
            + t
        \right).
    \end{equation}
\end{lemma}

\begin{lemma}[Concentration of $\norm{H_{[s]}}_1$] \label{lm:concentration_l1_Hs}
    With probability at least $1-2\delta$, we have
    \begin{equation}
        \abs{ \norm{\Hvecs}_1 - s \EE\left[Z | Z \geq t \right] }
        \leq
        \gencstpos \left( \sqrt{s}  + \sqrt{\log(1/\delta) s} \right).
    \end{equation}
\end{lemma}

\begin{lemma}[Concentration of $v_s$] \label{lm:concentration_vs}
    Assume $s < d/5$.
    For 
    $\delta \geq \exp(-s)$,
    with probability at least $1-6\delta$, we have
    \begin{align}
        \abs{
            \norm{v_s}_2^2 - s \left( \frac{2}{t^2} - \frac{10}{t^4} \right)
        } 
        &\leq s \left( \frac{\gencst}{t^6} +\Gencst{s,\delta} \right) ~\mathrm{and} \\
        \abs{
            \norm{v_s}_1 - s \left( \frac{1}{t} - \frac{2}{t^3} \right)
        }
        &\leq s \left( \frac{\gencst}{t^5} + \frac{\Gencst{s,\delta}}{t} \right) ~\mathrm{and} \\
        \abs{\innerprod{v_s}{H} - s}
        &\leq s \Gencst{s,\delta}
    \end{align}
    with
    $\Gencst{s,\delta} = \gencstpos \frac{t+ t\sqrt{\log(1/\delta)}}{\sqrt{s}}$.
\end{lemma}

The proposition follows as a consequence of this last lemma:
\begin{proof}[Proof of Proposition~ \ref{prop:concentration_gammas}]
    Let $c_1, c_2, \olt$ as in Lemma~\ref{lm:tlogds}, and assume $s \leq d/c_1$.
    In particular, $\log(d/s) \leq t^2 \leq 2\log(d/s)$.
    
    We apply 
    Lemma~\ref{lm:concentration_vs}
    with 
    $\delta = \exp \left( - 2\frac{s}{\log(d/s)^5} \right)$.
    Since $t^2 \leq 2 \log(d/s)$,
    this choice ensures that
    \mbox{$\frac{t\sqrt{\log(1/\delta)}}{\sqrt{s}} \leq 8/t^4$}.
    Moreover, the assumption that $d \leq \exp\left( \cd s^{1/5} \right)$ ensures that
    \mbox{$\frac{t}{\sqrt{s}} \leq \gencstpos/t^4$}.
    So we have
    \mbox{$\Gencst{s,\delta} \leq \gencstpos/t^4$}.
    
    
    The proposition follows by substituting the estimates of $\norm{v_s}_2^2$, $\norm{v_s}_1$, $\innerprod{v_s}{H}$ into 
    ${
        \gamma(\alphathr{s}) = \frac{\norm{H}_\infty}{\innerprod{v_s}{H}} v_s
    }$,
    and making the appropriate simplifications.
\end{proof}

\subsubsection{Proofs of the concentration lemmas}

\paragraph{Proof of Lemma~\ref{lm:concentration_Hs}: Concentration of $H_s$.}

By observing that
the random variable ${\max \{s ; H_s > t\} }$
is binomially distributed with parameters $d$ and $p = \PP(\abs{Z}>t)$,
\cite{li_2020} show the following upper and lower tail bounds for $H_s$.
\begin{claim} \label{claim:con_hs_t}
    Assume that $s < d/2$ and let $t$ be such that $2\Phic(t) = s/d$. Then for all $\eps>0$, we have the lower resp.\ upper tail bounds
    \begin{align}
        \PP(H_s \leq t - \eps) 
        &\leq \exp\left(
            -\gencstpos s \eps^2 \log(d/s) 
        \right) \\
        \mathrm{and}~~~
        \PP(H_s \geq t + \eps)
        &\leq \exp\left(
            -\gencstpos s \eps^2 \log(d/s)
            \exp\left( -2 \eps \sqrt{2\log(d/s) - \log\log(d/s) -\log(\frac{\pi}{2})}  - \eps^2 \right)
        \right).
    \end{align}
\end{claim}
\begin{proof}
    This follows straightforwardly from Lemma 2 of \cite{li_2020} and 
    from the estimate of $t$ in Lemma \ref{lm:tlogds}.
    %
\end{proof}
The lower tail bound is already sufficiently tight to show our high-probability lower bound on ${H_s-t}$.
However we remark that the upper tail bound is too loose; indeed it is only reasonable when $\eps$ is sufficiently small.
So to prove our high-probability upper bound, we instead start from the following one-sided concentration inequality from \cite{boucheron_2012}.
\begin{claim} \label{cr:uppertailbound_Hs}
    Assume 
    $d \geq 3$ and
    $s < d/2$,
    then 
    for all $z>0$,
    \begin{equation}
        \PP \left( H_{s} - \EE H_{s} \geq \gencstpos (\sqrt{z/s} + z/s) \right)
        \leq \exp(-z).
    \end{equation}
\end{claim}
\begin{proof}
    The proof follows from  the same argument as in Proposition 4.6 of \cite{boucheron_2012}.
\end{proof}

It only remains to bound the distance between $t$ and $\EE H_s$.
\begin{claim}
    Assume that $s < d/2$ and let $t$ be such that $2\Phic(t) = s/d$. Then
    \begin{equation}
        \abs{ \EE H_s - t } \leq \gencstpos \frac{1}{\sqrt{s}}.
    \end{equation}
\end{claim}
\begin{proof}
    According to Proposition 4.2 of \cite{boucheron_2012},
    \begin{equation}
        \Var(H_s) \leq
        \frac{1}{s \log 2}
        \frac{8}{
            \log\frac{2d}{s} - \log(1 + \frac{4}{s} \log\log \frac{2d}{s} )
        }
    \end{equation}
    so by Chebyshev's inequality,
    \begin{equation}
        \PP \left( \abs{H_s - \EE H_s} > \eps' \right)
        \leq \frac{\gencstpos'}{s} \frac{1}{(\eps')^2}.
    \end{equation}
    On the other hand, recall from Claim~\ref{claim:con_hs_t} that
    \begin{equation}
        \PP \left( \abs{H_s - t} > \eps \right)
        \leq
        2 \exp\left(
            -\gencstpos s \eps^2 \log(d/s)
            \exp\left( -2 \eps \sqrt{2\log(d/s) - \log\log(d/s) - \log(\frac{\pi}{2})}  - \eps^2 \right)
        \right)  
    \end{equation}
    One can check that there exist universal constants $\cone, \ctwo$ such that, by picking $\eps = \cone/\sqrt{s \log(d/s)}$ and $\eps' = \ctwo/\sqrt{s}$,
    the sum of the right-hand sides is less than $1$.
    
    Thus, with positive probability we have
    \begin{equation}
        \abs{\EE H_s - t} \leq \abs{H_s - \EE H_s} + \abs{H_s - t} \leq \frac{\cone / \sqrt{\log(2)} +\ctwo}{\sqrt{s}}.
    \end{equation}
\end{proof}

\paragraph{Proof of Lemma~\ref{lm:concentration_l2_Hs}: Concentration of $\norm{H_{[s]}}_2^2$.}

Let us first restate Proposition 2 of \cite{li_2020} in our notation.
We remark that their statement contained an additional $\log(d/s)$ factor due to a mistake in the proof. Correcting this mistake, we have that with probability at least $1-2\delta$,
\begin{equation}
    \abs{
        \frac{1}{\sqrt{s}} \norm{\Hvecs}_2 - \sqrt{\EE[Z^2|Z \geq t]}
    } \leq \gencstpos \frac{1}{\sqrt{s}} ( 1+\sqrt{\log(1/\delta)} ).
\end{equation}
Since for all $a,b,\eps>0$, 
$\abs{a-b} \leq \eps \implies \abs{a^2-b^2} \leq \eps (\eps + 2b)$,
this implies
\begin{equation}
    \abs{
        \frac{1}{s} \norm{\Hvecs}_2^2 - \EE[Z^2|Z \geq t]
    } \leq 
    \gencstpos \frac{1}{\sqrt{s}} (1+\sqrt{\log(1/\delta)})
    \left(
        \frac{1}{\sqrt{s}} (1+\sqrt{\log(1/\delta)})
        + \sqrt{\EE[Z^2|Z \geq t]}
    \right).
\end{equation}
Now 
$\EE[Z^2|Z \geq t] = 1+\frac{t^2}{h(t)} \leq \gencstpos t^2$
whenever $t \geq 1$, which is ensured by our assumption that $s/d = \Phic(t) \leq 0.2$.
So
\begin{equation}
    \abs{
        \frac{1}{s} \norm{\Hvecs}_2^2 - \EE[Z^2|Z \geq t]
    } \leq 
    \gencstpos \frac{1}{\sqrt{s}} (1+\sqrt{\log(1/\delta)})
    \left(
        \frac{1}{\sqrt{s}} (1+\sqrt{\log(1/\delta)})
        + t
    \right).
\end{equation}

\paragraph{Proof of Lemma~\ref{lm:concentration_l1_Hs}: Concentration of $\norm{H_{[s]}}_1$.}

We use exactly the same argument as in the proof of Proposition 2 of \cite{li_2020}.
Namely, start by decomposing
\begin{equation}
    \abs{
        \frac{1}{s} \norm{\Hvecs}_1 - \EE[Z | Z \geq t]
    }
    \leq
    \abs{
        \frac{1}{s} \norm{\Hvecs}_1 - \frac{1}{s} \EE \norm{\Hvecs}_1
    } + \abs{
        \frac{1}{s} \EE \norm{\Hvecs}_1 - \EE[Z | Z \geq t]
    }.
\end{equation}
For the first term, note that by rearrangement inequality, 
\mbox{$Z \mapsto \frac{1}{s}\sum_{i=1}^s \abs{Z_{(i)}}$} 
is $\frac{1}{\sqrt{s}}$-Lipschitz for the $\norm{\cdot}_2$ norm, 
where $(Z_{(1)},...,Z_{(d)})$ is the nondecreasing reordering of the absolute values of $Z$.%
\footnote{
    Proof: Denote $\tZ$ the reordering of $Z$ such that $\abs{\tZ_1}\geq...\geq\abs{\tZ_d}$ (but still with $\tZ$ signed).
    Then
    \begin{align}
        \abs{ \sum_{i=1}^s \abs{Z_{(i)}}-\abs{Y_{(i)}} } 
        = \abs{ \norm{\tZ_{[s]}}_1 - \norm{\tY_{[s]}}_1 }
        \leq \norm{\tZ_{[s]} - \tY_{[s]}}_1
        \leq \sqrt{s} \norm{\tZ_{[s]} - \tY_{[s]}}_2
        \leq \sqrt{s} \norm{Z - Y}_2
    \end{align}
    where the last inequality follows from the rearrangement inequality after taking squares and expanding.
}
So by concentration of Lipschitz-continuous functions of Gaussians,
\begin{equation}
    \PP \left(
        \abs{
            \frac{1}{s} \sum_{i=1}^s H_i - \EE \frac{1}{s} \sum_{i=1}^s H_i 
        }
        \geq \eps
    \right)
    \leq 2 \exp(-s \eps^2/2).
\end{equation}
For the second term, we can apply exactly the same arguments as in the proof of Proposition 2 of \cite{li_2020}, adapting equations (42) to (46), to obtain the bound
\begin{equation}
    \abs{
        \EE \frac{1}{s}\sum_{i=1}^s H_i - \EE[Z \vert Z \geq t]
    }
    \leq \gencstpos \EE \abs{H_{s+1} - t}
    \leq \gencstpos \frac{1}{\sqrt{s}}.
\end{equation}
In particular, we use the fact that 
$x \mapsto \EE [Z \vert Z \geq x]$
is a smooth function, which follows from its explicit expression given in Fact \ref{fact:cond_expectations}.

\paragraph{Proof of Lemma~\ref{lm:concentration_vs}: Concentration of $v_s$.}


For brevity of notation, let
$\Gencst{s,\delta} = \gencstpos \frac{t+ t\sqrt{\log(1/\delta)}}{\sqrt{s}}$.
Assume $\delta \geq e^{-s}$; 
in particular,
$\Gencst{s,\delta} \leq \gencstpos t$.
Collecting and simplifying the above results, so far we showed that
\begin{align}
    t \abs{H_s-t} 
    &\leq \Gencst{s,\delta} ~\mathrm{and}\\
    \abs{
        \frac{1}{s} \norm{\Hvecs}_2^2 - (1+\frac{t^2}{h(t)})
    } &\leq \Gencst{s,\delta} ~\mathrm{and} \\
    t \abs{
        \frac{1}{s} \norm{\Hvecs}_1 - \frac{t}{h(t)}
    } &\leq \Gencst{s,\delta}.
\end{align}

\begin{itemize}
    \item
    Substituting the deterministic estimates in the expression of $\norm{v_s}_2^2$ and carrying over the above concentration bounds, we obtain
    \begin{equation}
        \abs{
            \norm{v_s}_2^2 - \left(
            s(1+\frac{t^2}{h(t)})
            - 2 \frac{st}{h(t)} t
            + s t^2 \right)
        } \leq
        s \Gencst{s,\delta} (1+\Gencst{s,\delta})
    \end{equation}
    and the deterministic estimate can be simplified to
    \begin{equation}
        s \frac{(1+t^2) h(t) - t^2}{h(t)}
        = s \frac{
            \frac{2}{t^2} - \frac{12}{t^4} + O\left(\frac{1}{t^6}\right)
        }{
            1 - \frac{1}{t^2} + O\left(\frac{1}{t^4}\right)
        }
        = s \left(
            \frac{2}{t^2} - \frac{10}{t^4} + O\left(\frac{1}{t^6}\right)
        \right)
    \end{equation}
    (where the $O(\cdot)$ hides a universal constant).
    
    \item
    Likewise for $\norm{v_s}_1$ we get
    \begin{equation}
        \abs{
            \norm{v_s}_1 - \left( \frac{st}{h(t)} - st \right)
        } \leq 
        \frac{s}{t} \Gencst{s,\delta}
    \end{equation}
    and the deterministic estimate can be simplified to
    \begin{equation}
        st \left( \frac{1}{h(t)}-1 \right) 
        = s \left( \frac{1}{t^2} - \frac{2}{t^3} + O\left(\frac{1}{t^5}\right) \right).
    \end{equation}
    
    \item
    Likewise for $\innerprod{v_s}{H}$ we get
    \begin{equation}
        \abs{
            \innerprod{v_s}{H} - \left( s \left( 1 + \frac{t^2}{h(t)} \right) - \frac{st}{h(t)} t \right)
        } \leq 
        3 s \Gencst{s,\delta}
    \end{equation}
    and the deterministic estimate simplifies to $s$.
\end{itemize}

\subsubsection{Proof of Lemma~\ref{lm:tlogds}}
\label{apx:subsubsec:concentration_gamma:tlogs}

Using the upper bound 
\mbox{$\Phic(x) \leq \frac{\exp(-x^2/2)}{x \sqrt{2\pi}}$} from the first part of Fact~\ref{fact:estim_Phic},
it is straightforward to check that
$2\Phic(\olt) \leq s/d = 2\Phic(t)$,
and so $\olt \geq t$.

Let $\ult^2 = \olt^2 - \cone$ for some constant $\cone>0$ to be chosen.
Using the lower bound
\mbox{$\Phic(x) \geq \frac{\exp(-x^2/2)}{x \sqrt{2\pi}} \frac{x^2}{1+x^2}$},
one can check that $\cone$ can be chosen
such that
$2\Phic(\ult) \geq s/d = 2\Phic(t)$, and so $\ult \leq t$.
%

Going through the calculations reveals 
that $\cl \geq e^{2/\pi} \approx 2$ ensures $\olt^2 \leq 2\log(d/s)$,
that $2\Phic(\olt) \leq s/d$ is always true,
that 
$\cone = 
1 - \log(\frac{\pi}{2})
\approx 0.5$
ensures $\ult^2 \geq \log(d/s)$ for all $s, d$,
and that $\cl \geq e^{2.3415...} \approx 10.4$ ensures $2\Phic(\ult) \geq s/d$. 
This concludes the proof of the lemma. 

\begin{remark}
    Tighter bounds for $t$ can be derived while still only using the first-order estimate of $h(x)$ (the first part of Fact~\ref{fact:estim_Phic}).
    Namely, by similar straightforward calculations as above, one can check that
    there exist universal constants $\kappa, \alpha_1, \alpha_2>0$ such that, for all $s \leq d/\kappa$, $t$ is bounded as
    ${
        \ult \leq t \leq \olt
    }$
    where
    \begin{align}
        \ult^2 &= 2 \log(d/s) - \log\log(d/s) - \log(\pi)+ \frac{\log\log(d/s)}{2 \log(d/s)}  - \frac{\alpha_1}{\log(d/s)} \\
        \text{and}~~~
        \olt^2 &= 2 \log(d/s) - \log\log(d/s) - \log(\pi) + \frac{\log\log(d/s)}{2 \log(d/s)} + \frac{\alpha_2}{\log(d/s)}.
    \end{align}
\end{remark}

%% file: sections_proof/6_feasible_interval.tex
\subsection{Proofs of Lemmas~\ref{lm:interval} and \ref{lm:bounds}: Studying the feasible set of \texorpdfstring{\eqref{eq:phidp_reparam}}{A+'} and \texorpdfstring{\eqref{eq:phidm_reparam}}{A-'}}
\label{apx:subsec:feasible_interval}

\subsubsection{Proof of Lemma~\ref{lm:interval}
}
We give separate proofs for the statements 1-3 in the Lemma:
\newline\newline
\emph{First statement:  The mapping $\alpha \mapsto \norm{\gamma(\alpha)}_2^2$ is decreasing over $[1, \alphad]$ and increasing over $[\alphad, \alpham]$.} 

Using the notation of Section~\ref{apx:subsec:geometric_lemma}, the optimization problem defining $\gamma(\alpha)$ has Lagrangian 
    \begin{equation}
        L(w; \lambda, \mu, \nu) = \frac{1}{2} \norm{w}_2^2 
        - \lambda (\innerprod{w}{H} - \norm{H}_\infty)
        + \mu (\onevec^\top w - \alpha)
        - \nu^\top w,
    \end{equation}
    (up to the constant factor $\frac{1}{2}$ in the first term).
    By the envelope theorem, the marginal effect on the optimal value of increasing $\alpha$, is equal to the associated Lagrangian multiplier at optimum:
    $\frac{d \norm{\gamma(\alpha)}_2^2}{d\alpha} = -\mu$.
    A straightforward computation using the expression of $\mu$ from Lemma~\ref{lm:charact_gamma} reveals that
    $\mu>0$ for $1 < \alpha < \alphad$
    and $\mu<0$ for $\alphad < \alpha \leq \alpham$;
    hence the monotonicity of $\alpha \mapsto \norm{\gamma(\alpha)}_2^2$.
\newline \newline
\emph{Second statement: The mapping $\alpha \mapsto \norm{\gamma(\alpha)}_2^2$ is convex over $\left[1, \alpham \right]$, and $I$ is an interval.}

Recall that $\norm{\gamma(\alpha)}_2^2$ is given by
\begin{equation}
    \norm{\gamma(\alpha)}_2^2
    = \min_w  \norm{w}_2^2
    \subjto \begin{cases}
        \innerprod{w}{H} \geq \norm{H}_\infty \\
        \forall i, w_i \geq 0 \\
        \onevec^\top w = \norm{w}_1 = \alpha
    \end{cases}.
\end{equation}
Since $\alpha$ appears on the right-hand side of a linear constraint, it is straightforward to check directly that $\alpha \mapsto \norm{\gamma(\alpha)}_2^2$ is convex.
In detail: Let any $\alpha_0$ and $\alpha_1$, 
let $w_i = \gamma(\alpha_i)$ for $i \in \{1,2\}$, and $\alpha_t = (1-t) \alpha_0 + t \alpha_1$
and $w_t = (1-t) w_0 + t w_1$ for $t \in [0,1]$; then
$w_t$ is feasible for the optimization problem defining $\gamma(\alpha_t)$,
so $\norm{\gamma(\alpha_t)}_2^2
\leq \norm{w_t}_2^2
\leq t \norm{w_0}_2^2 + (1-t) \norm{w_1}_2^2$
by convexity.

$I$ is the $(B^2 \norm{H}_\infty^2)$-sublevel set of the function
$\alpha \mapsto (1-\cnoise) n \left( \sigma^2 \alpha^2 + B^2 \norm{\gamma(\alpha)}_2^2 \right)$,
which is convex, so $I$ is an interval.
\newline \newline
\emph{Third statement: The mapping $\alpha \mapsto \frac{\norm{\gamma(\alpha)}_2^2}{\alpha^2}$ is monotonically decreasing.} 

Note that for each $\alpha \in [1, \alpham]$, $\frac{\gamma(\alpha)}{\alpha}$ is the optimal solution of the optimization problem \begin{equation}
    \frac{\gamma(\alpha)}{\alpha} = \arg\min_w \norm{w}_2^2
    \subjto \begin{cases}
            \innerprod{w}{H} \geq \frac{\norm{H}_\infty}{\alpha} \\
            \forall i, w_i \geq 0 \\
            \onevec^\top w = \norm{w}_1 = 1
        \end{cases}.
\end{equation}
In particular, the constraint set is only increasing with $\alpha$, implying that $\alpha \mapsto \norm{\frac{\gamma(\alpha)}{\alpha}}_2^2$ is monotonically decreasing with $\alpha$.



\subsubsection{Proof of Lemma~\ref{lm:bounds}}
The upper bound for $\phidp$ immediately follows from Equation~\eqref{eq:phidp_reparam} and from the monotonicity of
$\alpha \mapsto \norm{\frac{\gamma(\alpha)}{\alpha}}_2^2$, which is the last statement of Lemma~\ref{lm:interval}.

For $\phidm$, there is an extra scale variable $0<\av \leq B$ which we first minimize out, similarly to the proof of Proposition~\ref{prop:main_norm} in Section~\ref{subsubsec:proof_bound_phidn}.
Starting from Equation~\eqref{eq:phidm_reparam}, first rewrite the constraint as
\begin{gather}
    \av^2 \norm{H}_\infty^2 \geq (1-\cnoise) n \left( \sigma^2 \alpha^2 + \av^2 \norm{\gamma(\alpha)}_2^2 \right) \\
    \iff
    \av^2 \geq 
     \underbrace{
         \frac{
            (1-\cnoise) n \sigma^2 \alpha^2
        }{
            \norm{H}_\infty^2 - (1-\cnoise) n \norm{\gamma(\alpha)}_2^2
        }
    }_{=: \fun(\alpha)^2}
    ~~~\text{and}~~~
    (1-\cnoise) n \norm{\gamma(\alpha)}_2^2 < \norm{H}_\infty^2.
\end{gather}

Then,
\begin{align}
    \phidm &=
    \min_{1 \leq \alpha \leq \alpham}
    \min_{\av}~
    \av^2 \frac{\norm{\gamma(\alpha)}_2^2}{\alpha^2}
    \subjto
    \fun(\alpha) \leq \av \leq B
    ~\text{and}~
    (1-\cnoise) n \norm{\gamma(\alpha)}_2^2 < \norm{H}_\infty^2 \\
    &= \min_{1 \leq \alpha \leq \alpham}
    \fun(\alpha)^2
    \frac{\norm{\gamma(\alpha)}_2^2}{\alpha^2}
    \subjto
    \fun(\alpha) \leq B
    ~\text{and}~
    (1-\cnoise) n \norm{\gamma(\alpha)}_2^2 < \norm{H}_\infty^2 \\
    &= \min_{1 \leq \alpha \leq \alpham}
    \fun(\alpha)^2
    \frac{\norm{\gamma(\alpha)}_2^2}{\alpha^2}
    \subjto
    \alpha \in I \\
    %
    &= \min_\alpha \frac{
        (1-\cnoise) n \sigma^2 \norm{\gamma(\alpha)}_2^2
    }{
        \norm{H}_\infty^2 - (1-\cnoise) n \norm{\gamma(\alpha)}_2^2
    }
    \subjto 
    \ulalpha_I \leq \alpha \leq \olalpha_I.
\end{align}
This objective is increasing in  $\norm{\gamma(\alpha)}_2^2$.
Now as shown in Lemma~\ref{lm:interval},
$\alpha \mapsto \norm{\gamma(\alpha)}_2^2$ is decreasing over $[1, \alphad]$.
Therefore, this objective is decreasing in $\alpha$ over 
$[\ulalpha_I, \alphad]$,
and
$\phidm$ is lower-bounded by its value at any $\olalpha$ such that
$\olalpha_I \leq \olalpha \leq \alphad$.
%

%% file: sections_proof/7_bound_phip_phim.tex
\subsection{Proof of Lemma~\ref{lm:boundsonalpha}: A tight admissible choice for \texorpdfstring{$\ulalpha$}{underlinealpha} and \texorpdfstring{$\olalpha$}{overlinealpha}}
\label{apx:subsec:bound_phip_phim}

Recall that we defined  
$ B(n,d) = c \sigma \sqrt{\norm{\wgt}_0} + M(n,d)$ where $\Brho(n,d)$ is given by Proposition~\ref{prop:main_norm} and  $c>0$ some universal constant.  For brevity of notation,  abbreviate $B(n,d) = B$. 
We seek integers $\us,\os$ such that $\alphathr{\us}$ and $\alphathr{\os}$ lie on the left, respectively on the right of the interval $I$, and such that $\alphathr{\os} \leq \alphad$, with high probability over the draws of $\Hgaussian$.

\subsubsection{Preliminaries}
We recall the notations
\mbox{$B = \Brho(n,d) + 2\norm{w^*}_1$} and
\begin{gather}
    I = \left\lbrace
        \alpha \in \left[ 1, \alpham \right]
        \middle\vert
        B^2 \norm{H}_\infty^2 \geq (1-\cnoise) n \left( \sigma^2 \alpha^2 + B^2 \norm{\gamma(\alpha)}_2^2 \right)
    \right\rbrace \\
    \text{and}~~~~~
    \forall \alpha \in [1, \alpham],~
    \fun(\alpha)^2 = \frac{
        (1-\cnoise) n \sigma^2 \alpha^2
    }{
        \norm{H}_\infty^2 - (1-\cnoise) n \norm{\gamma(\alpha)}_2^2
    }.
\end{gather}
Note that we have the equivalence
\begin{equation}
\label{eq:charactI_balpha}
    \alpha \in I \iff 
    (1-\cnoise) n \norm{\gamma(\alpha)}_2^2 < \norm{H}_\infty^2
    ~~~\text{and}~~~
    \fun(\alpha) \leq B.
\end{equation}

\paragraph{Reference point: $\alphathr{n} \in I$.}

We have by construction that, conditionally on the event where the inequalities in Equation~\ref{eq:concentration} hold for $s=n$, 
$\alphathr{n} \in I$ -- 
in particular this holds with probability at least $1-6\exp\left( -2\frac{n}{\log(d/n)^5} \right)$.
Indeed, let us take a closer look at the way we chose $\Brho(n,d)$,
from the proof of Proposition~\ref{prop:main_norm} (Section~\ref{subsubsec:proof_bound_phidn}).
We showed that, conditionally on that event,
$\alphathr{n}$ satisfies
\begin{equation}
    (1+\cnoise) n \norm{\gamma(\alphathr{n})}_2^2 < \norm{H}_\infty^2
    ~~~\text{and}~~~
    \ftild(\alphathr{n}) \leq \Brho(n,d)
\end{equation}
where
$
    \ftild(\alpha)^2
    = \frac{
        (1+\cnoise) n \sigma^2 \alpha^2
    }{
        \norm{H}_\infty^2 - (1+\cnoise) n \norm{\gamma(\alpha)}_2^2
    }
$.
Since
$\fun(\alphathr{n}) \leq \ftild(\alphathr{n})$
and $\Brho(n,d) \leq B$,
clearly $\alphathr{n}$ satisfies the condition~\eqref{eq:charactI_balpha}, \ie, $\alphathr{n} \in I$.

\paragraph{Summary of (in)equalities to be used in the proof.}
For ease of presentation, let us recall some assumptions or definitions that we will use throughout this proof.
For each integer $s$, $t_s \in \RR$ denotes the quantity such that $2\Phic(t_s) = s/d$, and $t_s^2 \asymp \log(d/s)$ by Lemma~\ref{lm:tlogds}.
By assumption, 
${
    \norm{\wgt}_0 \leq \cwone  \frac{\sigma^2 n}{\log(d/n)^{5}}
}$
for some universal constant $\cwone>0$.
By definition,
$\Brho(n,d)^2 =
\frac{\sigma^2 n}{t_n^2} 
\left(
    1 - \frac{2}{t_n^2} + \frac{\cbrho}{t_n^4}
\right)$
for some universal constant $\cbrho>0$.
In particular, this implies that
\begin{equation}
\label{eq:upperestim_B2}
    B^2 = \left( 
        \Brho(n,d) + c \sigma \sqrt{\norm{\wgt}_0}
    \right)^2 
    =
    \frac{\sigma^2 n}{t_n^2} 
    \left(
        1 - \frac{2}{t_n^2} + O\left(\frac{1}{t_n^4}\right)
    \right).
\end{equation}

\subsubsection{Finding \texorpdfstring{$\us$}{underlines} such that  \texorpdfstring{$\alphathr{\us} \leq \ualpha_I$}{alphaunderlines leq underlinealphaI}}


We want to find an $\us$ such that $\alphathr{\us}$ is on the left of the interval $I$, i.e., $\alphathr{\us} \leq \ualpha_I$.
Conditioning on the event where $\alpha_n \in I$,  it suffices  to have
$\alphathr{\us} < \alphathr{n}$ \ie\ $\us < n$,
and $\alphathr{\us} \not\in I$ \ie
\begin{align}
    &\frac{1}{B^2 \|H\|_\infty^2}
    (1-\cnoise) n 
    \left(
        \sigma^2 \norm{\gamma(\alphathr{\us})}_1^2
        + B^2 \norm{\gamma(\alphathr{\us})}_2^2) 
    \right) \\
    =&~ (1-\cnoise) n \frac{\sigma^2}{B^2} \frac{\norm{\gamma(\alphathr{\us})}_1^2}{\norm{H}_\infty^2}
    + (1-\cnoise) n \frac{\norm{\gamma(\alphathr{\us})}_2^2}{\norm{H}_\infty^2}
    > 1.
\end{align}

Instead of working directly with $\us$, it is more convenient to define $\us$ implicitly through a condition on $t_{\us}$.
Namely, we choose $\us$ such that $t_{\us}^2 \approx t_n^2 + \frac{\lambda}{t_n}$ for some constant $\lambda>0$. 
We now make this choice formal and show that $\us$ is very close to $n$; the following step will be to show that
this choice guarantees $\alpha_s \not\in I$.
\begin{claim} \label{claim:approxt_ulalpha_vgw}
    Assume $\cl n \leq d$.
    Let any fixed constant $0<\lambda\leq
    \sqrt{\log(\cl)}$, and
    let $\us$ be the largest integer $s$ such that $t_{s}^2 \geq t_n^2 + \frac{\lambda}{t_n}$.
    Then,
    \begin{gather}
        \us = n \exp\left( -\frac{\lambda}{2t_n} \right) \left( 1+O\left(\frac{1}{t_n^2}\right) \right) \\
        \text{and}~~~~~
        \abs{t_{\us}^2 - \left( t_n^2 + \frac{\lambda}{t_n} \right)}
        \leq O \left(
            \frac{1}{n}
        \right).
    \end{gather}
\end{claim}
The first equation quantifies the fact that this choice of $s$ is very close to $n$;
the second equation controls the error due to rounding (due to the fact that there is no integer $s$ such that $t_s^2 = t_n^2 + \frac{\lambda}{t_n}$ exactly).

\begin{proof}
   For concision, in this proof we write $s$ instead of $\us$. Denote $\olt_n^2 = t_n^2 + \frac{\lambda}{t_n}$.
    By definition, $t_s^2 \geq \olt_n^2 > t_{s+1}^2$.

    For the first part of the claim, we apply Fact~\ref{fact:estim_Phic} several times.
    Firstly,
    \begin{align}
        2\Phic(t_s) = \frac{s}{d}
        \leq \frac{2}{\sqrt{2\pi}} \frac{1}{t_s} e^{-t_s^2/2}
        \leq
        \frac{2}{\sqrt{2\pi}} \frac{1}{t_n} e^{-\olt_n^2/2}
        = 
        \frac{2}{\sqrt{2\pi}} \frac{1}{t_n} e^{-t_n^2/2} \exp \left( -\frac{\lambda}{2 t_n} \right).
    \end{align}
    Secondly,
    \begin{align}
        \frac{t_n^2}{1+t_n^2} \cdot
        \frac{2}{\sqrt{2\pi}} \frac{1}{t_n} e^{-t_n^2/2}
        &\leq 2 \Phic(t_n) \\
        \mathrm{and~hence}~~~~~
        \frac{2}{\sqrt{2\pi}} \frac{1}{t_n} e^{-t_n^2/2}
        &\leq \frac{n}{d} \left( 1 + \frac{1}{t_n^2} \right).
    \end{align}
    This proves the upper bound on $s$.
    The lower bound can be proved in a similar fashion, by applying Fact~\ref{fact:estim_Phic}
    to lower-bound $\frac{s+1}{d} = 2\Phic(t_{s+1})$,
    and again to upper-bound $\frac{1}{t_n^2} e^{-t_n^2/2}$.
    In this way we obtain
    \begin{equation}
        \frac{s+1}{d} \geq
        \frac{n}{d} \exp\left( -\frac{\lambda}{2 t_n} \right)
        \left( 1-\frac{1}{1+t_{s+1}^2} \right),
    \end{equation}
    and the bound $t_s^2 - t_{s+1}^2 \leq \frac{2}{s} \leq 2$ shown below implies that 
    $\frac{1}{1+t_{s+1}^2} = O \left(
        \frac{1}{t_n^2}
    \right)$.
    
    We now turn to the second part of the claim.
    By mean value theorem applied on $\Phic$, there exists $\xi \in [t_{s+1}, t_s]$ such that
    \begin{align}
        \frac{\Phic(t_{s+1}) - \Phic(t_s)}{t_{s+1}-t_s} 
        = \frac{1}{2d} \frac{1}{t_{s+1}-t_s}
        &= (\Phic)'(\xi)
        = -\frac{1}{\sqrt{2\pi}} e^{-\xi^2/2} \\ \mathrm{and}~~~~~
        0 < t_s - t_{s+1} 
        &= \frac{\sqrt{2\pi}}{2d} e^{\xi^2/2}
        \leq \frac{\sqrt{2\pi}}{2d} e^{t_s^2/2}.
    \end{align}
    Now, by Fact~\ref{fact:estim_Phic}, this can be further upper-bounded as
    \begin{align}
        \frac{\sqrt{2\pi}}{2d} e^{t_s^2/2}
        \leq \frac{1}{2d} \frac{1}{\Phic(t_s)} \frac{1}{t_s}
        = \frac{1}{s t_s}.
    \end{align}
    So we have the bound:
    \begin{equation}
        t_s^2 - \olt_n^2 
        \leq t_s^2 - t_{s+1}^2
        \leq \frac{t_s+t_{s+1}}{s t_s} 
        \leq \frac{2}{s}.
    \end{equation}
    (This completes the proof of the first part of the claim, for the lower bound.)
    We can conclude by substituting $s$ by its estimate from the first part of the claim,
    noting that $\frac{\lambda}{t_n}$ is uniformly bounded by assumption since
    $\lambda \leq \sqrt{\log(\cl)} \leq \sqrt{\log(d/n)} \leq t_n$
    by Lemma~\ref{lm:tlogds} (for an appropriate choice of $\cl$).
\end{proof}

We now show that we can choose the constant $\lambda>0$ such that
$\alphathr{\us} \not\in I$ with high probability. 

\begin{claim} \label{claim:ulalpha_notinI}
    The constants $\cn,\cl,\cd, \lambda>0$ can be chosen such that for any $n,d$ with $n \geq \cn$ and
    $\cl n \leq d \leq \exp(\cd n)$,
    \begin{equation}
        (1-\cnoise) n \frac{\sigma^2}{B^2} \frac{\norm{\gamma(\alphathr{\us})}_1^2}{\norm{H}_\infty^2}
        + (1-\cnoise) n \frac{\norm{\gamma(\alphathr{\us})}_2^2}{\norm{H}_\infty^2}
        > 1
    \end{equation}
    with probability at least 
    $1 - 12\exp \left( -\frac{n}{\log(d/n)^5} \right)$
    over the draws of $\Hgaussian$,
    where 
    $\us$ is defined  as in Claim~\ref{claim:approxt_ulalpha_vgw}.
\end{claim}

\begin{proof}
    We will repeatedly use the following inequalities, which follow from Lemma~\ref{lm:tlogds} and Claim~\ref{claim:approxt_ulalpha_vgw}:
    \begin{gather}
        t_n^2 = \log(d/n) + O(\log \log(d/n)) \\
        t_{\us}^2 = \log(d/n) + O(\log \log(d/n)) \\
        \frac{t_n^2}{t_{\us}^2} 
        = \frac{1}{1+ \frac{\lambda}{t_n^3} + O\left(\frac{1}{t_n^2n}\right)} 
        = 1 - \frac{\lambda}{t_n^3} + O\left(\frac{1}{t_n^2n}\right) + O\left(\frac{\lambda^2}{t_n^6}\right).
    \end{gather}
    
    Note that for appropriate choices of $\cn,\cl,\cd$, Equation~\eqref{eq:concentration} in Proposition~\ref{prop:concentration_gammas} holds simultaneously for $s=n$ and for $s=\us$ with probability at least
    $1 - 6\exp(-2\frac{n}{\log(d/n)^5}) - 6\exp(-2\frac{\us}{\log(d/\us)^5})
    \geq 1 - 12 \exp(-\frac{n}{\log(d/n)^5})$.
    We condition on this event throughout the remainder of the proof. 

    We begin with the first term, where we use the upper estimate of $B^2$ from Equation~\eqref{eq:upperestim_B2}:
    \begin{align}
        (1-\cnoise) n \frac{\sigma^2}{B^2} \frac{\norm{\gamma(\alphathr{\us})}_1^2}{\norm{H}_\infty^2}
        %
        &\geq 
        \frac{t_n^2}{1 - \frac{2}{t_n^2} + O\left(\frac{1}{t_n^4} \right)}
        \frac{1}{t_{\us}^2} \left( 1- \frac{4}{t_{\us}^2} + O\left(\frac{1}{t_{\us}^4}\right) \right) \\
        &\geq
        \left(1 - \frac{2}{t_n^2} + O\left(\frac{1}{t_n^4}\right)\right)
        \left( 1 - \frac{\lambda}{t_n^3} + O\left(\frac{1}{t_n^2 n} \right) + O\left(\frac{\lambda^2}{t_n^6} \right) \right) \\
        &= 1 - \frac{2}{t_n^2} 
        - \frac{\lambda}{t_n^3}  + O\left(\frac{1}{t_n^4}\right)
        + O\left(\frac{\lambda^2}{t_n^6}\right)
        +  O\left(\frac{\lambda}{t_n^2n}\right).
    \label{eq:proofphitmp1}
    \end{align}
    Next the  second term:
    \begin{align}
        (1-\cnoise) n \frac{\norm{\gamma(\alphathr{\us})}_2^2}{\norm{H}_\infty^2}
        &\geq n \frac{2}{\us t_{\us}^2}\left(1+ O\left(\frac{1}{t_{\us}^2}\right)\right) \\
        &\geq \frac{2}{t_{\us}^2} 
        \exp\left(\frac{\lambda}{2t_n}\right) 
        \left(1 + O\left( \frac{1}{t_n^2}\right) \right) \\
        &\geq \frac{2}{t_n^2} \left(1+ \frac{\lambda}{2t_n}  +  \frac{\lambda^2}{4t_n^2}+ O\left(\frac{\lambda^3}{t_n^3}\right)\right)\left(1 + O\left(\frac{1}{t_n^2}\right)\right).
    \label{eq:proofphitmp2}
    \end{align}
    Summing up \eqref{eq:proofphitmp1} and \eqref{eq:proofphitmp2}, we get:
    \begin{align}
        (1-\cnoise) n \frac{\sigma^2}{B^2} \frac{\norm{\gamma(\alphathr{\us})}_1^2}{\norm{H}_\infty^2}
        + (1-\cnoise) n \frac{\norm{\gamma(\alphathr{\us})}_2^2}{\norm{H}_\infty^2}
        \geq
        1 
        + \frac{\lambda^2}{2 t_n^4}
        + O\left(\frac{1}{t_n^4}\right)
        + O\left(\frac{\lambda^2}{t_n^6}\right)
        + O\left(\frac{\lambda}{t_n^2 n}\right). \label{eq:hereweneedwgt}
    \end{align}
    Clearly, we can choose the constants $\lambda, \cn,\cl,\cd>0$
    such that the right-hand side
    of the above equation
    is strictly greater than $1$ for any $n,d$ with
    $\cl n \leq d$, since $t_n^2 \asymp \log(d/n)$.
\end{proof}

\subsubsection{Finding \texorpdfstring{$\os$}{overlines} such that \texorpdfstring{$\oalpha_I \leq \alphathr{\os} \leq \alphad$}{overlinealphaI leq alphaoverlines}}

 We take the exact same approach to find $\os \geq n$ such that  $\alphathr{\os}$ is  on the right of the interval $I$, i.e., $\oalpha_I \leq \alphathr{\os}$. The derivations can be straightforwardly adapted, and we get the analogous results:
\begin{claim} \label{claim:approxt_olalpha_vgw}
    Assume $\cl n \leq d$.
    Let any fixed constant $0<\lambda \leq \sqrt{\log(\cl)}$, and
    %
    let $\os$ be the smallest integer $s$ such that $t_{s}^2 \leq t_n^2 - \frac{\lambda}{t_n}$.
    Then,
    \begin{align}
        \os &= n \exp\left( \frac{\lambda}{2t_n} \right) \left( 1+O\left(\frac{1}{t_n^2}\right) \right)
        & &\text{and} &
        \abs{t_{\os}^2 - \left( t_n^2 - \frac{\lambda}{t_n} \right)}
        &\leq O \left(
            \frac{1}{n}
        \right).
    \end{align}
\end{claim}
\begin{claim} \label{claim:olalpha_notinI}
The constants $\cn,\cl,\cd, \lambda>0$ can be chosen such that for any $n,d$ with $n \geq \cn$ and
    $\cl n \leq d \leq \exp(\cd n)$,
    \begin{equation}
        (1-\cnoise) n \frac{\sigma^2}{B^2} \frac{\norm{\gamma(\alphathr{\os})}_1^2}{\norm{H}_\infty^2}
        + (1-\cnoise) n \frac{\norm{\gamma(\alphathr{\os})}_2^2}{\norm{H}_\infty^2}
        > 1
    \end{equation}
    with probability at least 
    $1 - 12\exp \left( -\frac{n}{\log(d/n)^5} \right)$
    over the draws of $\Hgaussian$, 
    where $\os$ is defined as in Claim~\ref{claim:approxt_olalpha_vgw}.
\end{claim}


It remains to check that this choice of $\os$ satisfies $\alphathr{\os} < \alphad$.
But this is clearly the case, because $\os<d$ (by the first part of Claim~\ref{claim:approxt_olalpha_vgw}, for appropriate choices of $\cl$) and $\alphathr{d} < \alphad$ by definition.
This concludes the proof of Lemma~\ref{lm:boundsonalpha}.